\documentclass[12pt,leqno, titlepage]{article}

\usepackage{epsfig}
\usepackage{amssymb} \textwidth=6.5in \hoffset=-0.75in
\usepackage[all]{xy}
\usepackage[active]{srcltx}
\usepackage{euscript}
\usepackage{graphicx}
\usepackage{psfrag}
\usepackage{color}
\usepackage{amsmath}

\newtheorem{theorem}{Theorem}[section] 
\newtheorem{lemma}{Lemma}[section] \baselineskip 1.5em
\newtheorem{corollary}{Corollary}[section]
\newtheorem{proposition}{Proposition}[section]
\newtheorem{definition}{Definition}[section]
\newtheorem{remark}{Remark}[section]

\newcommand{\norm}[1]{\left\Vert#1\right\Vert}
\newcommand{\R}{\mathbb R}

\newcommand{\nablash}{\nabla{\kern -.75 em
     \raise 1.5 true pt\hbox{{\bf/}}}\kern +.1 em}
\newcommand{\Deltash}{\Delta{\kern -.69 em
     \raise .2 true pt\hbox{{\bf/}}}\kern +.1 em}
\newcommand{\Rslash}{R{\kern -.60 em
     \raise 1.5 true pt\hbox{{\bf/}}}\kern +.1 em}

\newcommand{\Hb}{\bar H}

\newcommand{\D}{\partial}

\newcommand{\id}{\operatorname{id}}
\newcommand{\Real}{\operatorname{Re}}

\title{{\large Anisotropic Einstein data with isotropic nonnegative scalar curvature}}

\author{Bernold Fiedler*\\
Juliette Hell*\\ 
Brian Smith*\\[2mm]
*Institut f\"ur Mathematik\\
Freie Universit\"at Berlin\\
Arnimallee 3, D--14195 Berlin, GERMANY}

\begin{document}

\maketitle

\setlength{\parindent}{0cm}

\section{Introduction}

The problem of constructing 3-manifolds of prescribed non-negative scalar curvature is 
important for the initial value formulation 
of the Einstein equations in general relativity in the asymptotically flat case.  
For some standard references see \cite{bartnikmax, cafimamu, cantorini}. For a more recent survey article 
involving the constraint equations see \cite{BartIsen}.
In the maximal ``gauge'', valid initial 
data consist of an asymptotically flat Riemannian manifold $(M,g)$ together with a symmetric, trace free tensorfield $k_{ij}$, vectorfield $J$,
and non-negative function $\rho$.  Here $k_{ij}$ is interpreted as the second fundamental form of $(M,g)$ as embedded in spacetime, 
$J$ is the mass current density, and $\rho$ is the local mass density; the non-negativity of the latter follows from the dominant energy condition.    
In addition, the data $(M,g,k,J,\rho)$
must satisfy the Einstein constraint equations, which in the maximal ``gauge'' appear as 
 \begin{eqnarray} \label{einsteinconstraint}
       R(g)-|k|^2_g&=&16\pi\rho\\
     \nabla\cdot k&=&-8\pi J,
\end{eqnarray}
where $R(g)$ denotes the scalar curvature of the metric $g$ and $\nabla$ refers to the Levi-Civita connection with respect to $g$.
Thus, the construction of such initial data involves the construction of a manifold $(M,g)$ of non-negative scalar curvature.  
In the time symmetric case $k=0$ the scalar curvature is proportional to the mass density $R(g)=16\pi\rho$, and the problem reduces to a 
prescribed non-negative scalar curvature problem.      

\begin{figure}[h]
\includegraphics[width=0.9\textwidth]{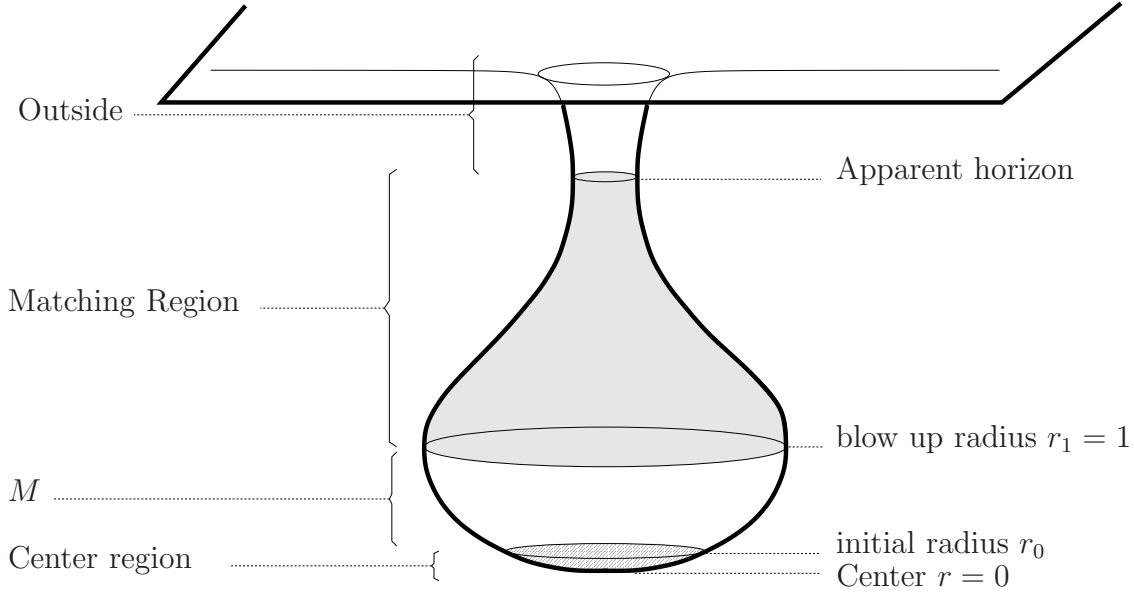}
\caption{ \label{vase}  {\em 2d-caricature of a black hole initial data for the Einstein equations. 
The vertical variable is the radius $r$. Each horizontal section $\{r= \rm{constant} \}$ 
in this figure is a circle, standing in fact for a 2-sphere $\Sigma=S^2$.  
The region outside the apparent horizon (and above it in this picture) was constructed in \cite{Smith1, Smith2}. 
The present paper constructs the region called $M$ between an initial radius $r_0>0$ and a blow up radius $r_1=1$ 
whose section $\{ r_1\}\times S^2$ is a critical point of the area functional. The matching and the center regions (shaded) 
are still to be constructed.}}
\end{figure}

The  prevalent method for solving the constraints has been the conformal method. This requires at the outset a metric $\tilde g$ 
which is conformally equivalent to a metric of  non-negative scalar curvature. Such a metric $\tilde{g}$  
is only guaranteed  to exist if  the $L^{3/2}$ norm of the negative part of the 
scalar curvature $R(\tilde g)$ is small~\cite{cabri}. That
is, the scalar curvature of $\tilde g$ must be almost non-negative already.  Thus, what is needed is a non-conformal method for constructing 
general metrics of prescribed non-negative scalar curvature.  Solving the parabolic scalar curvature equation (\ref{eq:parscal}) below
 provides just such a method. It is our main goal, in the present paper, to construct many solutions of the parabolic curvature 
equation by means of equivariant bifurcation theory and symmetry breaking. For previous ad-hoc constructions see \cite{S}.  
Figure \ref{vase} sketches a 2-dimensional caricature of the 3-dimensional space initial data for a black hole $M$ 
foliated over a radial variable $r$ by 2-dimensional spheres $\Sigma=S^2$, 
which are only horizontal circle in Figure \ref{vase}. 
\psfrag{mR}{Matching Region}
\psfrag{B}{Outside}
\psfrag{M}{$M$}
\psfrag{C}{Center region}
\psfrag{c}{Center $r=0$}
\psfrag{i}{initial radius $r_0$}
\psfrag{bu}{blow up radius $r_1=1$}
\psfrag{E}{Apparent horizon} 

Given a foliation of the manifold $M$, with foliating function $r$, the parabolic scalar curvature equation relates the 
scalar curvature to the transversal 
metric component $u=|\nabla_g r|_g^{-1}$.  We will only be concerned with regions upon which $r$ is non-critical so
that we may assume without loss of generality that the manifold takes the form $M=[r_0,r_1]\times\Sigma$, where $\Sigma$ 
is a regular 2-manifold - in the present work we take $\Sigma=S^2$, but the equation is true for any $\Sigma$.   
After a foliation preserving diffeomorphism, any metric $g$ on $M$ can be expressed as
\begin{equation} \label{riemmetric}
   g=u^2dr^2+r^2\omega,
\end{equation}
where $\omega$ is a family of metrics on $\Sigma$ extended to $M$ so that $\omega(\D_r,\cdot)=0$.
The parabolic scalar curvature equation can then be expressed as  
\begin{equation} \label{eq:parscal}
     \Hb r\D_r u =u^2\Delta_{\omega}u+Bu-\left(\kappa-\frac{1}{2}r^2R\right)u^3,
\end{equation}
where 
\begin{equation} \label{Hbar}
      \Hb=2+r\partial_r\omega_{ab}\omega^{ab},
\end{equation}
$\Delta_{\omega}$ is the Laplace-Beltrami operator of $\omega$, $\kappa$ 
is the Gauss curvature of $\omega$, and $B$ is also 
determined entirely by $\omega$; specifically 
$B=r\D_r\Hb-\Hb+\frac{1}{2}(\Hb^2+|\omega+r\D_r \omega|^2_{\omega})$. 
Equation~(\ref{eq:parscal}) was first derived in the quasi-spherical case by  Bartnik~\cite{bartnik93}. 
It was fully generalized to the 3-dimensional case by Smith and Weinstein~\cite{SW}.  
It was later extended to $n$-dimensions by Shi and Tam~\cite{ST}. 

Concerning the type of  equation (\ref{eq:parscal}), note that the mean curvature of the foliation is given by
$H=\Hb/ru$.  Thus, assuming that the family $\omega$ is such that the mean curvature of the foliation is positive, 
equation~(\ref{eq:parscal}) is in fact parabolic, with the radial variable $r$ playing the role of the time variable. 
Accordingly, we shall often speak of $r$ and increasing functions of $r$ as time variables in analogy with the heat equation, 
porous medium equation, etc. For convenience, we also call the  boundary value problem with prescribed data at $\{r_0\}\times\Sigma$
an initial value problem. 

Assuming compactness of $\Sigma$, one uses this equation to construct a metric $g$ of prescribed scalar curvature $R$ in the following way:  
choose a smooth family of metrics such that $\Hb>0$, as well as ``initial'' data $u_0\in C^{\infty}(\Sigma)$, $0<u_0$.  
Solve the parabolic scalar curvature equation ~(\ref{eq:parscal}) with initial data $u(r_0)=u_0$. By standard parabolic regularity theory 
this is always possible for a small enough interval $[r_0,r_0+\varepsilon)$.  

Of course,  even in the case of long time existence, 
the resulting manifold will not be complete since we have a boundary at $r=r_0$.  Assuming $\Sigma=S^2$, it is possible 
to produce solutions starting at 
$r=0$~\cite{bartnik93}, but we shall not worry about this in the present work. For small $r\in [0,r_0]$ it is simpler for the metric to have been constructed
by other means.  For instance,  one can first construct the metric in Gaussian normal coordinates for $r<r_0$, and then use the 
parabolic scalar curvature equation to extend to $r\geq r_0$.  For more interesting topologies, one will need to start with an inner minimal surface,
which will then be joined to a similarly constructed manifold on the other side.    
For more on this see~\cite{Smith1},~\cite{Smith2}.

Concerning global existence, it suffices that
the operator 
\begin{equation} \label{opT}
    T=\Delta_{\omega}-\left(\kappa-\frac{r^2R}{2}\right)
\end{equation}
be non-positive~\cite{Smith1}. Since we are primarily interested in asymptotically flat manifolds whose ends have the topological structure of
$[r_0, +\infty)\times S ^2$, we choose, in fact, $\Sigma=S^2$.  Then by the Gauss-Bonnet theorem $\int_{S^2}\kappa=4\pi$, and 
we are allowed to choose $\omega$ such that $\kappa>0$.  Positivity of the operator $-T$ can be ensured by choosing $R<2r^2\kappa$.

Nonetheless, there are many situations 
in which this condition  is not satisfied that do, in fact, lead to blow-up at finite $r$.  If the mass density 
$\rho=R/16\pi$ is very large so that, for instance
$\kappa-\frac{r^2R}{2}\leq -c<0$, one can easily see from the maximum principle that blow-up will occur at finite time (radius) $r_1$.  This is, of 
course, only blow-up of a metric component, but it is geometrically significant with our choice of foliation: the mean curvature of the 
foliation approaches 
$0$ at any blow-up point, and  parabolicity is no longer satisfied at the blow-up ``time'' $r=r_1$. This will be true for any 
bounded ``time'' variable that we might use.      

It may, however, happen that the constructed metric itself retains sufficient regularity at the blow-up ``time''  $r=r_1$ so as to allow an
extension beyond this radius by other means.  The simplest example occurs when $\omega$ is just the standard metric $\mathring{\omega}$ of $S^2$ and $r^2R\equiv 4$. 
Then the parabolic scalar curvature equation becomes 
\begin{equation} \label{eq:triv}
     2r\D_r u = u^2\Delta u +u^3+u,
\end{equation}
where we have written $\Delta=\Delta_{\mathring\omega}$.
On $[r_0,1)\times S^2$, for instance, this admits the trivial blow-up solution 
\begin{equation} \label{trivialu}
u=\left(1/r-1\right)^{-1/2}, 
\end{equation}
with corresponding metric 
\begin{equation} \label{trivialblowupmetric}
   g=\left(1/r-1\right)^{-1}dr^2+r^2\mathring\omega
\end{equation}
Obviously the metric blows up at $r=1$.  In terms of  the geodesic distance function 
\begin{equation} \label{geodistance}
     s=\int_{r_0}^r \left( 1/\tilde r-1\right)^{-1/2} d\tilde r,
\end{equation}
however the metric takes the form $g=ds^2+r^2(s)\mathring\omega$.  
It is easily checked that $r$ is $C^{\infty}$ smooth  in $s$, and thus $g$ is
$C^{\infty}$ smooth up to and including the boundary at $r=1$.

A natural next line of investigation, which we continue to follow in this work, 
is to ask if there exist anisotropic solutions $u$ of  (\ref{eq:parscal}) in the case 
$\omega=\mathring\omega$  with the same property. Here anisotropic means, in other words, that  $u=u(r,p)$ 
is not only a function of $r$ alone but also depends on the angular variables $p\in \Sigma=S^2$. Such nonhomogeneous 
solutions are not invariant under the action 
of the full symmetry group $O(3)$ any more: they exhibit anisotropy.  More specifically, 
we would like to address the existence of self-similar blow-up solutions in the form 
\begin{equation} \label{rescalingu}
     u=\left(1/r-1\right)^{-1/2}\nu,
\end{equation}
for a function $\nu=\nu(r,p)$, which is bounded above and below by positive constants, and is in addition $C^{\infty}$ 
on the interior  $(r_0, 1)\times S^2$ and $C^0$ at the boundary $ \{r_0, 1\}\times S^2$.  Using the function $s$ defined in the 
previous paragraph (which is now 
only proportional to the geodesic distance), the corresponding metric $g$ can be expressed as
\begin{equation} \label{zemetric}
     g= \nu^2 ds^2 +r^2(s)\mathring\omega,
\end{equation}
and the aforementioned regularity of $\nu$ is exactly reflected in the metric $g$.  In this work we shall not always obtain regularity at the 
boundary $r=1$ beyond $C^0$, 
but it should in any case be noted that in order to do so it is more appropriate to use a function such as $s$ for the foliating 
function since the area radius
variable $r$ is degenerate at $r=1$.  In fact, in the present case we can say more:  the mean curvature $H$ of the 
outer boundary  surface $r=1$ vanishes 
identically; i.e. it is a minimal surface.  Even in the case that further differentiability at $r=1$ is lacking, 
the stability functional - alias the second variation of the area functional- is
defined and shows this surface to have at least one unstable direction - more on this in Section 7.

The approach described in the preceding paragraph was initiated in~\cite{S}.  In that work it is shown that if the function $\nu$ 
is  bounded away from $0$, then it is also uniformly bounded above, and thus the final constructed metric is 
$C^{\infty}([r_0,1)\times S^2)\cap  C^{0}([r_0,1]\times S^2)$.  Furthermore, under suitable conditions on 
$r^2R$, non-linear stability of the trivial 
blow-up was obtained. In particular this showed the existence of at least some non-trivial solutions of (\ref{eq:parscal})
whose blow up behavior is only asymptotically given by the trivial solution (\ref{trivialu}).
However, this could certainly not be described as a plethora of solutions. Part of the aim of the present 
work is to find more blow-up solutions in this category by a more systematic bifurcation analysis.  


Define $\lambda=\left(r^2R-2\right)$.  In this work we take $\lambda$ to be an adjustable constant, i.e. a bifurcation parameter.  
Substituting  (\ref{rescalingu}) $u=(1/r-1)^{-1/2}\nu$   in equation (\ref{eq:parscal})
we obtain the equation for $\nu$:
\begin{equation} \label{eq:nu,r}
     2(1-r)\partial_r \nu = \nu^2\Delta \nu+\frac{\lambda}{2} \nu^3 - \nu.  
\end{equation}
We may eliminate the breakdown of parabolicity at $r=1$ by defining a new variable 
$t$  such that $r=1-\exp(-2t)$.  The surfaces $r\rightarrow 1$ are described by $ t \to\infty$ in the self-similarly rescaled equation. 
\begin{equation} \label{eq:nueq}
     \partial_{t} \nu = \nu^2\Delta \nu+\frac{\lambda}{2} \nu^3 - \nu.   
\end{equation}
A few observations are in order: $C^0$ Regularity of the metric $g$ at $r=1$ corresponds 
to the existence of a globally bounded solution $\nu$, $0<\mu\leq \nu\leq M<\infty$, 
converging to an equilibrium $\nu_*\in C^{\infty}(S^2)$ of equation (\ref{eq:nueq}) as $t \to\infty$.  
The equilibrium  $\nu_*$ is a solution of the equation 
\begin{equation} \label{eq:nuequil}
    \Delta \nu_*+\frac{\lambda}{2} \nu_* - \frac{1}{\nu_*}=0.   
\end{equation}
It will be convenient  to rescale, $\nu \mapsto \nu/\sqrt{\lambda/2},\, t\mapsto t\lambda/2$, 
so that the equilibrium equation  (\ref{eq:nuequil}) becomes 
\begin{equation}  \label{eq:nuequil2}
        \Delta\nu_*+\frac{\lambda}{2}\left(\nu_* - \frac{1}{\nu_*}\right)=0,
\end{equation}
and the dynamical equation becomes 
\begin{equation}  \label{eq:nudyn}
       \D_{t}\nu=\nu^2\left(\Delta\nu+\frac{\lambda}{2}\left(\nu - \frac{1}{\nu}\right)\right),
\end{equation}

In order to apply analytic semigroup theory and equivariant bifurcation theory, it is more convenient to translate 
the trivial equilibrium solution $\nu _*\equiv 1$ to the origin by defining  $v=\nu-1$:
\begin{equation}\label{eq:v1.*}
     \D_{t} v=(v+1)^2 \left( \Delta v+\lambda  f(v) \right). 
\end{equation}
Here
\begin{equation} \label{func}
 f(v)=v-\frac{1}{2}v^2/(1+v).
 \end{equation}
We often refer to equation~(\ref{eq:v1.*}) as the rescaled equation in the following. The equilibria $v_*$ of the rescaled equation verify
\begin{equation}\label{eq:vequil}
      \Delta v_*+\lambda f(v_*)=0.
\end{equation}
Note that such an equilibrium solution corresponds to self-similar blow-up 
\begin{equation}\label{ch:tildenutou}
     u(r,p)=\left(1/r-1\right)^{-1/2}\nu_*(p), \qquad p\in S^2.   
\end{equation}

Rephrased in the new notation, the work~\cite{S} only treats the case $\lambda=1$ explicitly: 
trivial self-similar  blow-up corresponding to $v_*\equiv 0$, and a local strong stable manifold of solutions which converge to this.    
In the present work we 
vary $\lambda$ and use $O(3)$ equivariant bifurcation theory           
to produce branches of solutions. Indeed, as $\lambda$ crosses the eigenvalues $\lambda_{\ell}=\ell(\ell+1)$ of $-\Delta_{S^2}$, 
the linearized left hand side of the equation~(\ref{eq:vequil}) becomes linearly degenerate, and  $O(3)$ 
equivariant bifurcation theory produces branches of solutions along the appropriate isotropy subspaces.
At each branch point we obtain one or several local symmetry breaking families of solutions $v_*$ of the equilibrium equation.  
Each such solution $v_*$ yields a self-similar metric 
metric  
\begin{equation} \label{metselfsim}
    g = \frac{\lambda}{2}\frac{(1+v_*)^2}{\left(1/r-1\right)}dr^2+r^2\mathring\omega,
\end{equation}
which can be seen to be $C^{\infty}$ smooth up to and including the boundary $r=1$ in the radially geodesic coordinates $s$.  These anisotropic metrics do not posses full $O(3)$ symmetry, 
but rather the symmetries of the isotropy subspaces of the eigenspace  of $\Delta_{S^2}$, 
which were encountered at the branch point from which it originated.  

Based on~\cite{L1}, we adapt a 
strong stable manifold theorem to apply at each of these new equilibria.  
This allows us to construct immortal and eternal solutions  $v$ to the translated rescaled equation (\ref{eq:v1.*}). Immortal, i.e. 
defined and bounded uniformly for $t \rightarrow +\infty$, alias $r\rightarrow 1$, are solutions $v$ in the strong 
stable manifolds of equilibria $v_*$. 
These converge to $v_*$ and correspond to asymptotically self-similar solutions $u$, as indicated in (\ref{ch:tildenutou}, \ref{metselfsim}) above. 
Eternal solutions $v$ exist for all real $t$ with uniform bounds, and are heteroclinic between different equilibria: 
$v\rightarrow v_{\pm}$ for $t \rightarrow \pm \infty$. The metric interpretation $u$ of ancient solutions $v$       which exist and are globally bounded for $t \rightarrow - \infty$, alias $r\rightarrow -\infty$, and which 
constitute the unstable manifolds of equilibria $v$ of (\ref{eq:v1.*}), will be discussed in Section 7.

Equilibria $v_*$ and solutions $v$ asymptotic to them in fact exhaust the dynamic possibilities of uniformly bounded eternal solutions to (\ref{eq:v1.*}). 
Indeed the standard energy functional 
\begin{equation}
E(v):= \int_{S^2} \left(  \frac12 \vert \nabla v \vert^2 - \lambda F(v)  \right) dp,
\end{equation}
where $F$ is a primitive function of the nonlinearity $f$, i.e. $\D_v F=f$. The functional $F$ is a Lyapunov function for the semilinear variant 
\begin{equation}  \label{semvar}
\D_{t} v = \Delta v + \lambda f(v)
\end{equation}
of (\ref{eq:v1.*}). Indeed (\ref{semvar}) can be interpreted as the $L^2$ gradient semi-flow of the Lyapunov functional $E(v)$ because 
\begin{equation}
\frac{d}{dt} E(v) = - \int_{S^2} v_{t}^2 dp.
\end{equation}
In particular, uniformly bounded solutions tend to equilibrium in any time direction. Similarly, the rescaled curvature equation(\ref{eq:v1.*}) can be 
interpreted as the gradient semi-flow of $E(v)$ with respect to a slightly adapted $L^2$-metric which depends explicitly on $v$. Specifically
\begin{equation} \label{func2}
\frac{d}{dt} E(v) = - \int_{S^2} (1+v)^{-1} v_{t}^2 dp.
\end{equation}
for solutions $v$ of (\ref{eq:v1.*}) with uniformly positive weight $ (1+v)^{-1}$. In particular uniformly bounded solutions still 
tend to equilibrium, in any time direction. See \cite{S} where accumulation on several equilibria is excluded using \cite{lsimon}.  


Our main results, Proposition \ref{prop7.2}, Theorem \ref{th7.3} and Table \ref{tableheteroclinics} below, can be summarized as follows. 
Equation  \ref{eq:parscal}, where $\omega$ is the standard metric on $S^2$ independently of $r$ and $\lambda=r^2R-2$.
We obtain  blow-up solutions of the form 
\begin{equation} \label{solmainth}
 u(r,p)  = \left(\frac{\lambda}{2}\left( 1/r -1\right)\right)^{-\frac12}\left( v(-\frac{1}{2}\log(1-r), p) + 1\right), \qquad p\in S^2,
 \end{equation}
where $v$ is an anisotropic function on the 2-sphere whose remaining symmetry is described by the isotropy groups of Table \ref{tabletypes}. 
Figure \ref{fig7.1} shows how those functions bifurcate from the fully isotropic constant function on $S^2$  as the parameter $\lambda$ varies. 
Table \ref{tableheteroclinics} describes in each case the heteroclinic connections between bifurcating equilibria. 
The corresponding anisotropic metrics on $M$ are regular up to the blow up radius and exhibit there a minimal surface. 


In sections 2, 3, we summarize some necessary background  from equivariant bifurcation theory.
In Section 4, we check that one of the main technical ingredients - the strong stable manifold -  is provided in our quasilinear case, 
as well as the semi-group framework. In Section 5, we specify the isotropy subgroups which play a role in the symmetry 
breaking bifurcations of equation (\ref{eq:v1.*}). The rest of the paper puts all these elements together to obtain the aforementioned result;
see Sections 6 and 7.  

This work was supported by the Deutsche Forschungsgemeinschaft, SFB 647 "Space--Time--Matter".


\section{Symmetry and equivariance}\label{sym}

We briefly recall some symmetry terminology for equivariant
dynamics. In the subsequent Section \ref{branching}, we formulate the equivariant
branching lemma which is an elementary but very useful variant of a
classical bifurcation theorem of Crandall and Rabinowitz \cite{CrRab}. 
ce
 For the convenience of the reader the precise role of
equivariance will be emphasized.
More generally see \cite{GolStew, Vdbwequiv, CLbook}
for a
background on symmetry and Lyapunov-Schmidt reduction, and \cite{Carr, Vdbw, L1,  koch} for center manifolds in semilinear settings. 

Let $X,Y$ be Banach spaces with bounded linear group actions
$\varrho^X, \varrho^Y$ of the same group $\Gamma$. In other words
$\varrho^X,\varrho^Y$ are group homomorphisms from $\Gamma$ to the
invertible bounded linear operators on $X,Y$, respectively:

\setcounter{equation}{0}

\begin{equation}
\varrho^X(\gamma_1 \cdot \gamma_2) = \varrho^X(\gamma_1)\cdot
\varrho^X(\gamma_2)\label{eq2.1} 
\end{equation}
for all $\gamma_1,\gamma_2 \in \Gamma$, and similarly for $Y$. We call
$\varrho^X,\varrho^Y$ {\it representations} of $\Gamma$ on
$X,Y$. Frequently we will use the abbreviated notation $\gamma u:=
\varrho^X(\gamma) u$, for $\gamma \in \Gamma$ and $u \in X$.
We call a representation $\varrho^X$ strongly continuous if $(\gamma, v)\mapsto \varrho (\gamma)v$ is continuous. For finite-dimensional $X$ this is equivalent to continuity of homomorphism $\varrho^X:\Gamma \rightarrow GL(X)$ from $\Gamma$ to the general linear group on $X$.  

 A map 
\begin{equation}
F: X \to Y \label{eq2.2}
\end{equation}
is called $\Gamma$-{\it equivariant} under $\varrho^X,\varrho^Y$, if
$F(\varrho^X(\gamma)v) = \varrho^Y(\gamma)F(v)$ holds for all
$\gamma\in\Gamma, v \in X$. To simplify notation we simply write this
requirement as
\begin{equation}
F(\gamma v) = \gamma F(v) \label{eq2.3}
\end{equation}
without too much ambiguity. In the important special case $X \subseteq
Y$ where $X$ is a subspace of $Y$, albeit with a possibly stronger
norm, it is particularly convenient to forget the distinction between
$\varrho^X$ and $\varrho^Y$, when $\varrho^X$ simply restricts
$\varrho^Y$ to $X$.

The intuitive notion of ``symmetry'' can be made precise in two
slightly different ways. Given $v \in X$ and a representation
$\varrho^X$ of $\Gamma$ we call
\begin{equation}
\Gamma_v:= \{\gamma \in \Gamma \text{\ } \vert \text{\ } \gamma v = v\}, \label{eq2.4}
\end{equation}
i.e. the group of elements $\gamma$ which fix $v$,  the {\it isotropy}
of $v$ or the {\it stabilizer} of $v$. For example $\Gamma_0=\Gamma:$
trivially $v=0$ possesses full isotropy. Note that
\begin{equation}
\Gamma_{\gamma v} = \gamma \Gamma_v \gamma^{-1}; \label{eq2.5}
\end{equation}
i.e. elements of the same group orbit possess conjugate isotropy.

Conversely, given any subgroup $K$ of $\Gamma$, we can consider the
closed subspace
\begin{equation}
{\rm Fix}_X(K):= X^K:= \{v\in X\text{\ } \vert \text{\ }\gamma v = v \ \ {\rm for \ all} \ \
\gamma \in K\}\label{eq2.6}
\end{equation}
of $K$-{\it fixed vectors}, which we call the {\it fix space} of
$K$. For example $K \le \Gamma_v$ is a subgroup of the isotropy $\Gamma_v$, for
any $v \in {\rm Fix}(K)$.

Let $F:X \to Y$ be $\Gamma$-equivariant, as above, and fix any
subgroup $K \le \Gamma$. Then $F$ restricts to
\begin{equation}
F: X^K \to Y^K, \label{eq2.7}
\end{equation}
i.e. $F$ maps $K$-fixed vectors $v \in X$ to $K$-fixed vectors $F(v)
\in Y$. Indeed, let $\gamma \in K, v \in X^K$. Then $\gamma F(v) =
F(\gamma v) = F(v)$ and hence $F(v) \in Y^K$. In this general way the fix spaces
${\rm Fix} (K)$ capture the idea of an ``ansatz'' to respect some symmetry
property $K$, quite efficiently and generally.

The isotropy subgroups of a representation $\varrho^X$ of $\Gamma$ on
$X$ form a lattice as follows. Vertices are the conjugacy classes of
isotropy subgroups. Containment of conjugacy classes of subgroups defines the lattice
structure.

A representation $\varrho^V$ of $\Gamma$ on a real or complex Banach
space $V$ is called {\it irreducible} if $\varrho^V$ does not restrict
to a representation on any nontrivial closed subspace $\{0\} \neq
\tilde{V} \neq V$. In other words $\varrho ^V(\tilde{V})$ is never
contained in $\tilde{V}$, for any such $\tilde{V}$. By Schur's Lemma,
the only $\Gamma$-equivariant {\it complex linear} maps on a
finite-dimensional complex irreducible representation are the complex
multiples of identity. Therefore irreducible representations play a
central role in spectral analysis.

For example consider the Sobolev space $X=H^2(S^2)$, $Y= L^2(S^2)$ and
the selfadjoint Laplace-Beltrami operator $-\Delta$ on $Y$ with domain
${\cal D}(\Delta) = X$. Then $- \Delta : X \to Y$ is equivariant under
the action of the compact group of orthogonal matrices $\gamma \in
\Gamma:= O(3)$ on $v \in Y$ given by
\begin{equation}
(\gamma v)(x) = (\varrho^Y(\gamma) v) (x) := v(\gamma^{-1}x). \label{eq2.8}
\end{equation} 
The representation restricts to $X=H^2 \subseteq L^2$. The spectrum
${\rm spec}(- \Delta)$ is given by real eigenvalues $0 < \lambda_0 <
\lambda_1 <\ldots$ with
\begin{equation}
\lambda_\ell = \ell(\ell+1).\label{eq2.9}
\end{equation}
The associated eigenspaces $V_\ell$ of dimension $2\ell+1$ enumerate
the irreducible continuous representations of the special orthogonal group $SO(3)$, for
$\ell=0,1,2,\ldots$. They are spanned, for each $\ell$, by the real
and imaginary parts of the spherical harmonics
\begin{equation}
Y_{\ell m}(\vartheta,\varphi)= P^m_\ell(\cos \vartheta)\exp(im \vartheta)\label{eq2.10a}
\end{equation}
with $m=0,\ldots,\ell$. The real Legendre functions $P^m_\ell$ 
are essentially $m$-th derivatives of Legendre polynomials $P_\ell(x)$
with a prefactor $(1-x^2)^{m/2}$. The spherical harmonics $Y_{\ell m}$ are normalized to
form a complete orthonormal basis of $Y=L^2(S^2)$, in polar or zenith
angle coordinates $0 \le \vartheta \le \pi$ and azimuth $0 \le \varphi
\le 2\pi$. Because the Legendre polynomials $P_\ell(x)$ are even/odd
for even/odd $\ell$, the reflection $\gamma = -id \notin SO(3)$ acts as multiplication by $(-1)^\ell$ on $V_\ell$.

As a first example we claim
\begin{equation}
{\rm Fix}_{V_{\ell}}(SO(2)) = \ \ {\rm span} \ \ Y_{\ell 0}. \label{eq2.10b}
\end{equation}
Here elements of $SO(2)$, and the rotations around the
polar axis $\vartheta=0$ of $S^2$, shift the azimuth $\varphi$ but
keep the polar angle $\vartheta$ fixed. By Fourier expansion with
respect to $\varphi$ this proves (\ref{eq2.10b}). We call functions
$v(\vartheta,\varphi) = v(\vartheta)$ in Fix$(SO(2))$ {\it
 axisymmetric}.

Note how the precise isotropy $\Gamma_v \le SO(3)$ of axisymmetric
functions $v$ depends on $\ell$. Let $\mathbb{Z}^c_2:= \{\pm \mathrm{id} \}$
denote the center subgroup of $\Gamma=O(3)= SO(3) \cup (-SO(3))$. Then
\begin{equation}
\Gamma_v = \left\{
\begin{array}{lcl}
O(2) \oplus \mathbb{Z}^c_2 & {\rm for} & \ell > 0 \  \ {\rm even}\\[2mm]
\quad O(2)^- & {\rm for} & \ell > 0 \ \ {\rm odd}
\end{array}
\right.\label{eq2.10c}
\end{equation}
for non zero $v = Y_{\ell 0}$. The group $O(2)^-$
 refers to an isotropy
subgroup such that 
$\Gamma_v$ is isomorphic to $O(2)\leq SO(3)$ but
$\Gamma_v\cap SO(3) = SO(2)$ is a subgroup of
index 2. Specifically $O(2)^-$ are the planar rotations $SO(2)$ in the
equatorial plane $\vartheta = \pi/2$ around the polar axis, and the
in-plane reflections $\varphi \mapsto \varphi_0-\varphi$
via reflections at perpendicular planes. The group $O(2) \le SO(3)$,
in contrast, achieves reflections in the equatorial plane by 180$^o$
rotations around in-plane axes with $\vartheta=\pi/2$. By $O(2) \oplus
\mathbb{Z}^c_2$ we denote the direct product of $O(2)$ with $\pm$id,
i.e. $O(2) \oplus \mathbb{Z}^c_2 = O(2) \cup (-O(2))$.

In Figure \ref{fig2.1} we sketch the lattices of (conjugacy classes
of) isotropy subgroups $\Gamma_v$ for the representation (\ref{eq2.8})
of $\Gamma = 0(3)$ on $V_\ell, \ \ell=0,\ldots,4$. In parentheses we
indicate $\dim_{\mathbb{R}}$Fix$(\Gamma_v)$ of $\Gamma_v$. A complete
classification had been initiated  by \cite{IhGol}(1984); see
also the detailed and corrected accounts in the monumental work of
Lauterbach \cite{LautHabil}  as contained in
Chossat Lauterbach Melbourne \cite{CLM} (1991); see also \cite{L89}.

A {\it maximal isotropy subgroup} $\Gamma_v$ is a strict subgroup
$\Gamma_v \neq \Gamma = O(3)$ which is maximal in the lattice with that
property. In the cases $\ell = 1,\ldots,4$ of figure 2.1 these
coincide with the isotropy subgroups $\Gamma_v$ of real
one-dimensional fix spaces, 
\begin{equation}
\dim {\rm Fix}\ \ (\Gamma_v) = 1.\label{eq2.11}
\end{equation}
Obviously condition (\ref{eq2.11}) implies 
that the isotropy subgroup is maximal: one-dimensional linear representations $\varrho$ of 
groups admit only actions $\varrho(\gamma)=\pm 1$. Conversely, however,  
it is not true that fix spaces of maximal isotropy subgroups are one-dimensional. 
The counterexample with lowest $\ell$ is $\ell=12$,  the octaedral group $\Gamma_v =O\oplus\ \mathbb{Z}_2^c$, and $\dim {\rm Fix}(\Gamma_v)=2$.  

For further details on the nomenclature of subgroups of $O(3)$ we refer to \cite{CLM, CLbook,  CL}. 
We only describe the maximal isotropy subgroups in Figure \ref{fig2.1} here. We have already 
explained the isotropies $O(2)^-$ and $O(2)\oplus\mathbb{Z}_2^c$ of axisymmetric $v\in {\rm span} Y_{\ell 0 }$ 
in (\ref{eq2.10c}). Since $-\mathrm{id} \notin SO(3)$ acts as $(-1)^{\ell}\mathrm{id}$ on $V_{\ell}$, any element of $V_{\ell}$ 
possesses isotropy at least
$\mathbb{Z}^c_2 = \{\pm \mathrm{id}\}$, if $\ell$ is even. Therefore any closed isotropy
subgroup $\Gamma_v$ of $\Gamma = O(3)$ 
takes the form $\Gamma_v = K \oplus \mathbb{Z}^c_2= K \cup (-K)$ with
$K \le SO(3)$. The subgroups of $SO(3)$ are $SO(2),O(2),T,O,I$,  the
discrete cyclic subgroups $\mathbb{Z}_n$ of $SO(2)$, and the dihedral
subgroups $D_n$ of $SO(2)$. Here $D_n$ is the subgroup of planar
rotations and reflections in $O(2) \le SO(3)$ which fix a regular
$n$-gon. Similarly $T,O,I \le SO(3)$ fix the regular tetrahedron,
octahedron, and icosahedron and are isomorphic to the alternating
symmetric group $A_4$, to the symmetric group $S_4$, and to the alternating symmetric group $A_5$,
respectively.

The ``octahedral'' group $O^- \le O(3)$ is the group of orthogonal
rotations {\it and reflections} which fix a regular tetrahedron: this time
the subgroup of rotations is the tetrahedral group $T = O^- \cap
SO(3)$ of index 2 in $O^-$. Similarly the dihedral group $D^d_6 \le
O(3)$, isomorphic to $D_6 \le O(2)$, has index 2 subgroup $D^d_6 \cap
SO(3) = D_3 \le O(2)$.

We summarize the results of this section, following \cite{LautHabil, L89, CLM, CLbook}.
\setcounter{figure}{0} 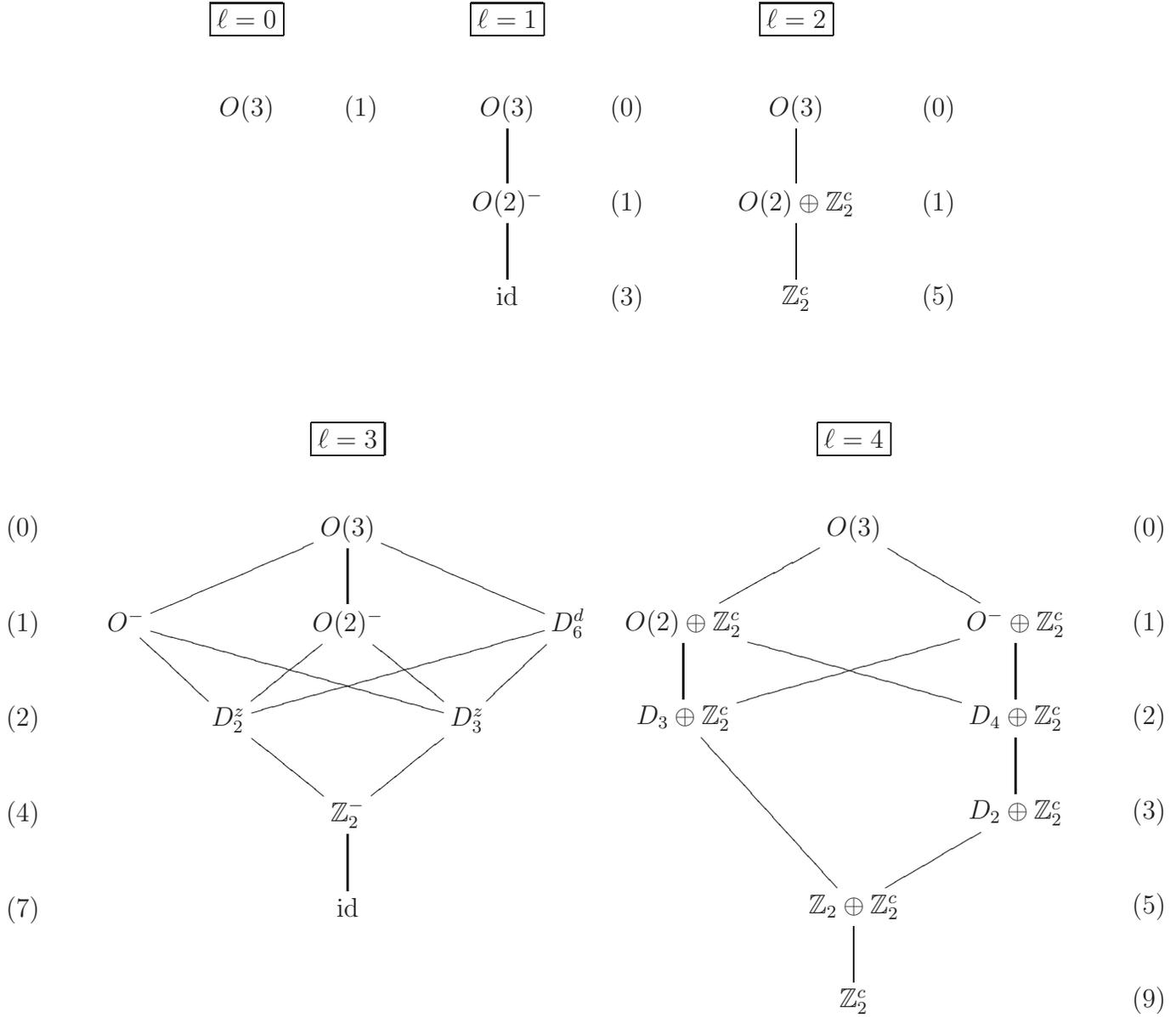
\begin{figure}[!t]
\centerline{
\begin{tabular}{c}
\begin{tabular}{ccccccc}
\xymatrix{
  *+[F]{\ell=0}& \\
 O(3) &  (1)
}
& \qquad &\qquad &
\xymatrix{
*+[F]{\ell=1} &\\
O(3) \ar@{-}[d] &(0)\\
O(2)^-\ar@{-}[d] & (1)\\
\rm{id} & (3)
}
&\qquad &\qquad &
\xymatrix{
*+[F]{\ell=2} & \\
O(3)  \ar@{-}[d]& (0) \\
O(2)\oplus \mathbb{Z}_2^c  \ar@{-}[d] & (1)\\
 \mathbb{Z}_2^c & (5)
}
\end{tabular}
\\
\\
\\
\\
\begin{tabular}{cc}
\xymatrix{
 &&&*+[F]{\ell=3} && \\
(0)&&& O(3)  \ar@{-}[dll]  \ar@{-}[d] \ar@{-}[drr]& &\\
(1) &O^-  \ar@{-}[dr]  \ar@{-}[drrr]& &O(2)^- \ar@{-}[dl]  \ar@{-}[dr]& &D_6^d\ar@{-}[dl]  \ar@{-}[dlll]\\
(2)& & D_2^z \ar@{-}[dr]  & & D_3^z\ar@{-}[dl]& \\
(4)&&&\mathbb{Z}_2^- \ar@{-}[d] &&\\
(7)&&& \rm{id} &&
}
& 
\xymatrix{
&*+[F]{\ell=4} &&\\
 & O(3) \ar@{-}[dl] \ar@{-}[dr] & & (0)\\
O(2)\oplus \mathbb{Z}_2^c \ar@{-}[drr]  \ar@{-}[d]  && O^-\oplus \mathbb{Z}_2^c \ar@{-}[dll]  \ar@{-}[d] & (1)\\
D_3\oplus \mathbb{Z}_2^c \ar@{-}[ddr]&&  D_4 \oplus \mathbb{Z}_2^c\ar@{-}[d]  & (2) \\
&& D_2\oplus \mathbb{Z}_2^c \ar@{-}[dl]& (3) \\
 & \mathbb{Z}_2 \oplus \mathbb{Z}_2^c\ar@{-}[d] & & (5) \\
 & \mathbb{Z}_2^c && (9)
}
\end{tabular}
\end{tabular}
}
\caption{\em{The lattices of conjugacy classes of isotropy subgroups
   of the natural irreducible representation of $O(3)$ on the spaces
   $V_\ell$ of spherical harmonics, $\ell =0,\ldots,4$. See \cite{CL}
   }}\label{fig2.1} 
\end{figure}
\begin{lemma} (Lauterbach et al.\cite{LautHabil, CLM, CLbook})\\
\label{lm2.1}
Consider the irreducible representations $(\ref{eq2.8})$ of $\Gamma =
O(3)$ on the spherical harmonics subspaces $V_\ell$; see
$(\ref{eq2.10a})$.

Then the lattices of conjugacy classes of isotropy subgroups
$\Gamma_v$ on $V_\ell$ and the real dimensions of the fix
spaces Fix$(\Gamma_v) = V^{\Gamma_v}_\ell$ are given in Figure
$\ref{fig2.1}$. The one-dimensional fix spaces of maximal isotropy
subgroups are given in Table \ref{tablegenerators}, for axisymmetric solutions and all
$\ell$, and for the remaining cases in $\ell = 3,4$.
\end{lemma}
\begin{table}[!h]
\begin{tabular}{||c||c|c|c|c|c||}
\hline \hline
&&&&&\\
$\Gamma_v$ & $O(2)\oplus \mathbb{Z}^c_2$ &  $O(2)^-$  &  $ D_6^d$ & $O^-$ &  $O\oplus \mathbb{Z}_2^c$  \\
&&&&&\\
\hline 
&&&&&\\
$\ell >0 $&even& odd&3&3&4 \\
&&&&&\\
\hline 
&&&&&\\
$\mathbf{e} \in$ Fix$(\Gamma_v)$&$Y_{\ell 0 }$&$Y_{\ell 0 }$&Re$Y_{33 }$&Re$Y_{32 }$& $Y_{40}+ \sqrt{5/14}Y_{42}$\\
&&&&&\\
\hline \hline
\end{tabular}\\
\caption{ Explicit generators $\mathbf{e}$ of the one-dimensional fix spaces
$\rm{Fix}(\Gamma_v)$ of  isotropy subgroups $\Gamma_v$ in certain
spherical harmonics representations $V_\ell$ of $\Gamma = O(3)$. \label{tablegenerators} }
\end{table}
The graphs over $S^2$ of the  generators $\mathbf{e}$ f the one-dimensional fix spaces
Fix$(\Gamma_v)$ of  isotropy subgroups $\Gamma_v$ listed in Table \ref{tablegenerators} are represented in Figure \ref{shapes}. 
\begin{figure}[h]
\begin{center}
\includegraphics[width=0.4\textwidth]{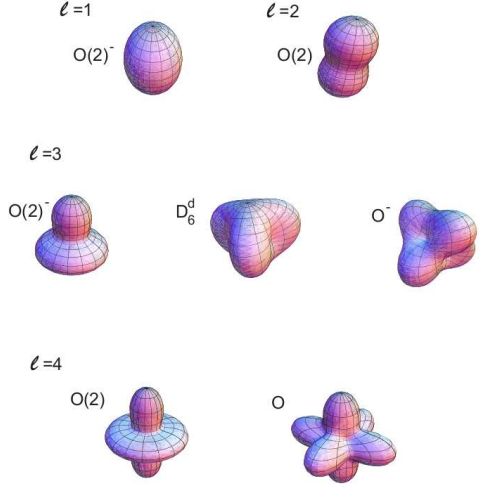}
\end{center}
\caption{\label{shapes}  {\em The graphs of the generators $\mathbf{e}$ of Table \ref{tablegenerators} over the sphere $S^2$}}
\end{figure}

\section{Equivariant branching}
\label{branching}
\setcounter{equation}{0}

The equivariant branching Lemma \ref{lm3.1} below was developed by Vanderbauwhede \cite{Vdbwequiv}
 into a simple, but strikingly effective tool of
equivariant bifurcation analysis; see also \cite{Cic}. As we will see,
it is a straightforward adaptation of the classical bifurcation
theorem by Crandall and Rabinowitz on bifurcation of zeros of maps \cite{CrRab}:

We consider $C^k$-maps, $k\geq 2$
\begin{equation}
 \label{eq3.1}
 \begin{array}{lll}
G: &&\mathbb{R}\times X \to Y  \\
&& (\lambda,v) \mapsto G(\lambda,v)
\end{array}
\end{equation}
with Banach spaces $X,Y$ and a trivial branch
\begin{equation}
G(\lambda,0) = 0\label{eq3.2}
\end{equation}
of zeros. Moreover we assume $\Gamma$-equivariance, i.e.
\begin{equation}
G(\lambda,\gamma v) = \gamma G(\lambda,v)\label{eq3.3}
\end{equation}
for all $\gamma$ in the group $\Gamma$, all $\lambda \in \mathbb{R}$ and $v
\in X$. In particular the observations of Section  \ref{sym} apply to $F=G(\lambda, v)$; see (\ref{eq2.3}). For any subgroup $K \le \Gamma$ let 
\begin{equation}
G^K: \mathbb{R} \times X^K \to Y^K\label{eq3.4}
\end{equation}
denote the restriction of $G$ to $\mathbb{R} \times X^K$, where
$X^K,Y^K$ are the $K$-fixed closed subspaces of $X,Y$ introduced in
(\ref{eq2.6}), (\ref{eq2.7}).

\begin{lemma}(Vanderbauwhede\cite{Vdbwequiv})\\
\label{lm3.1}
{~ }
In the setting $(\ref{eq3.1})$--$(\ref{eq3.4})$ assume there exists a
subgroup $K\le \Gamma$ and $\lambda_0 \in \mathbb{R}$ such that the
linearization $L^K(\lambda) = D_x G^K(\lambda,0)$ is Fredholm of 
index zero. Moreover assume
\begin{equation}
\ker L^K(\lambda_0) = \ {\rm span} \ \{\mathbf{e}\}\label{eq3.5}
\end{equation} 
is real one-dimensional, and the transversality condition
\begin{equation}
\partial_\lambda L^K(\lambda_0) \mathbf{e} \notin  \ {\rm Range} \ L^K(\lambda_0)\label{eq3.6}
\end{equation}
holds.

Then the nontrivial solutions of $G^K(\lambda,v) = 0$ near $(\lambda_0,0)$
form a bifurcating $C^{k-1}$  curve
\begin{equation}
s\mapsto (\lambda(s),v(s)) \in \mathbb{R} \times X^K \label{eq3.7}
\end{equation}
which satisfies
\begin{equation}
\lambda(0) = \lambda_0, \quad u(0) = 0, \quad v'(0) = \mathbf{e}. \label{eq3.8}
\end{equation}
In particular the bifurcating curve intersects the trivial branch $(\lambda,0)$
only at $(\lambda_0,0)$, where $s=0$. The isotropy $\Gamma_{v(s)}$ is
at least $K$, along the branch.
\end{lemma}

{\bf Proof:}
{~}

Statements (\ref{eq3.7}), (\ref{eq3.8}) 
follow directly from \cite{CrRab}.
The intersection claim is a consequence of $v'(0)=\mathbf{e}\neq 0$.
 The isotropy claim follows from $v(s)
\in X^K$; see the example following definition (\ref{eq2.6}) of
$X^K$. This proves the lemma.
\hfill$\bowtie$

\section{Center and strong stable manifolds}
\label{mnf}
\setcounter{equation}{0}

We shall now obtain invariant manifold theorems for Equation~(\ref{eq:v1.*})  which rewrites  as
\begin{equation} \label{v1.*in4}
       \partial_{t} v=(v+1)^2 \left( \Delta v+\lambda  f(v) \right)
\end{equation}
at an arbitrary equilibrium $v_*$.  We regard this as an abstract dynamical system in the Banach space $X=C^{\beta}(S^2), \beta\in (0,1)$.  
and follow the approach of~\cite{L1} closely. \\
Notation: $C^{\delta}$, where $\delta\in \mathbb{R}$, $\delta=[\delta]+\{\delta\}$, $[\delta]\in \mathbb{N}$ entire part of $\delta$, $\{\delta\}\in [0,1)$ its decimal part,  denotes the space of $[\delta]$-times differentiable functions whose $[\delta]$-th derivative is $\{\delta\}$-H\"{o}lder continuous. \\
We take  $X=C^{\beta}$ and 
$D=C^{2+\beta}(S^2)$ as the domain of the right hand side.  Let $\tilde v = v-v_*$.  We must rewrite the equation in terms of $v_*$ 
and $\tilde v$.  
This is not difficult since the right hand side is just a polynomial in $v,v_*,\Delta v,\Delta v_*$.  We shall not carry this out 
in full detail here, but shall simply note that the equation takes the form 
\begin{equation}   \label{eq:dyneq}
    \tilde v_t=A\tilde v+Q(\tilde v),  
\end{equation}
where $A$ is just the linearization of the right hand side at $v_*$, and $Q$ is the nonlinear part, which satisfies $Q(0)=0,Q'(0)=0$.  
The operator $A$ is easily computed to be 
\begin{equation}
      A= (1+v_*)^2\left(\Delta + \lambda\left(1+\frac{1}{(1+v_*)^2}\right)\right).
\end{equation}

To reduce the technical requirements of the present paper, at this stage, we observe that the linearization of 
(\ref{eq:v1.*}) or (\ref{eq:dyneq})
is Morse equivalent to the semilinear parabolic PDE 
\begin{equation} \label{eq:semilinear}
    \D_tv=\Delta v+ f(v),
\end{equation}
again with $f(v)=2v-v^2/(1+v)$ as in~(\ref{func}).  By Morse equivalent we mean that the total algebraic multiplicity 
$i(v_*)$ of the positive spectrum $\sigma_+$ of the linearization $L_*$ at any equilibrium $v_*$ of~(\ref{eq:semilinear}) coincides 
with that of the linearization $A$ of ~(\ref{eq:v1.*}). 
\begin{lemma}
Equilibria of $v_*$ of~(\ref{eq:v1.*}),~(\ref{eq:semilinear}) are Morse equivalent.  
\end{lemma}
{\bf Proof:}
{~}
   
 By the $L^2$ self-adjointness of $L_*=\Delta+\lambda f'(v_*)$ we notice that all eigenvalues $\mu$ of the linearization $L_*$ 
and of the Linearization 
$A=(1+v_*)^2L_*$ are real.  Indeed, $A w=\mu w$ in $L^2(S^2)$ implies 
\begin{equation}
 \mu\int_{S^2}(1+v_*)^{-2}w\bar w = \mu\left<(1+v_*)^{-1}w,w\right>=\left<L_*w, w\right>
=\left<w,L_*w\right>=\left<w,L_*w\right>=\left<w, L_* w\right>\in\R
\end{equation}
for the standard Hermite $L^2$ scalar product $\left<\cdot,\cdot\right>$ on $S^2$.  Hence $\mu\in\R$. Moreover, $\mu=0$ is 
an eigenvalue of $L_*$ if and only if $\mu=0$ is an eigenvalue of $A$ with equal multiplicity.  The homotopy of operators
\begin{equation}
    \tilde L(\tau)=(1+\tau v_*)^2L_*
\end{equation}
with $0\leq\tau\leq 1$ shows the same property.  Therefore the linearizations $L_*$ and $A$ are Morse equivalent and the lemma
is proved. 

Here it is important that the algebraic and geometric multiplicities coincide.  We prove this indirectly: suppose the algebraic
multiplicity of $\mu$ is strictly larger.  Then there exists a nonzero eigenvector $w_1$ with $\mu w_1=A w_1$ and $w_2$ such that
$w_1=(\mu-A)w_2$.  Abbreviating $m=1/(1+v_*)^2>0$, this implies 
\begin{align}
   0&=\left(m\mu-L_*\right)w_1 \label{m&mu}\\
   mw_1&=\left(m\mu-L_*\right)w_2.
\end{align}
Because $L_*$ and $\mu m$ are self-adjoint this implies 
\begin{equation}   \label{anotherdamnequation}
 \left<w_1,mw_1\right>=\left<w_1,(m\mu-L_*)w_2\right>=\left<(m\mu-L_*)w_1,w_2\right>=0.
\end{equation}
This contradiction to positivity of $mw_1^2$ proves that the geometric and algebraic multiplicities of $L$ coincide.  
\hfill$\bowtie$ 

The little lemma above enables us to derive Morse indices $i(v_*)$, as well as the multiplicity of any purely imaginary spectrum from 
considerations of only the semilinear parabolic equation (\ref{eq:semilinear}) instead of the quasilinear problem~(\ref{eq:v1.*}).  For 
semilinear semigroup settings, spectra and their associated center, center unstable and strong 
stable invariant manifolds have been studies very extensively in the existing literature.  See for instance~\cite{henry}, \cite{Vdbw},  \cite{brufie}, and the many references there. In particular, we will 
invoke center manifolds
in our local 
bifurcation analysis of (\ref{eq:v1.*}) in Section 5 and 6, to establish normal hyperbolicity and determine the Morse indices of nontrivial equilibria 
$v_*\neq 0$ of reduced isotropy, which bifurcate from the trivial equilibrium $v_*=0$ at $\lambda =\lambda_l=l(l+1)$.  

Once these Morse indices of normally hyperbolic equilibria $\Gamma v_*\neq 0$ have been determined, however, we need to study the
remaining solutions $v(t)$ of the original equation~(\ref{eq:v1.*}) which converge to $v_*$ for $t\to\infty$.  These solutions $v(t)$ 
form strong stable manifolds of $v_*$, i.e. manifolds characterized by an exponential decay with bounded exponentially weighted norm
\begin{equation}\label{eq:superfluous}
        \Vert v-v_* \Vert_{\varepsilon}   :=    \sup_{t\geq 0}\norm{v(t)-v_*}e^{\varepsilon t} 
\end{equation}
for small enough $\varepsilon$.  See~\cite{koch}, \cite{Mielke} for pertinent results in a quasilinear setting.   
Next, we adapt such a result in the setting of 
\cite{L1}.

As described previously, we work with equation (\ref{eq:dyneq})  
which is the result of a self-similar rescaling at the ODE blow up rate, 
and centering on an equilibrium $v^*$ of the rescaled equation. The ambient Banach space 
is chosen to be $X=C^{\beta} (S^2)$, $\beta \in (0,1)$. The domain $D$ of the linear operator 
$A$ is defined as $X\supset D := C^{2+\beta}(S^2)$. Our approach consists in applying Theorem 9.1.7 of 
\cite{L1}, and for that we introduce some basic interpolation spaces between $D$ and $X$, and 
state what they turn out to be in the case of the equation (\ref{eq:dyneq}) we are studying.

\begin{definition} \label{interp}
For $\alpha\in (0,1)$,  we define  
\[
    D_A(\alpha,\infty)=\left\{v\in X: t\mapsto\norm{t^{1-\alpha}Ae^{tA}v}_X\in L^{\infty}(0,1) \right\}
\]
and 
\[
     \norm{v}_{D_A(\alpha,\infty)}=\norm{v}_X+[v]_{D_A(\alpha,\infty)},
\]
where $[v]_{D_A(\alpha,\infty)}$ denotes the $L^{\infty}$ norm in time over $[0,1]$ 
of the defining function $\norm{t^{1-\alpha}Ae^{tA}v}_X$.
Furthermore, we recursively define 
\[
    D_A(1+\alpha,\infty) = \left\{v\in D: Av\in D_A(\alpha,\infty)\right\}, 
\]
with 
\[
   \norm{v}_{D_A(1+\alpha,\infty)}=\norm{v}+\norm{Av}_{D_A(\alpha,\infty)}.
\]
\end{definition}
\begin{remark}
The interpolation spaces defined in  \ref{interp} are Banach spaces if the operator $A$ is sectorial. In our settings, 
\begin{eqnarray} \label{ambientanddomain}
X&=&C^{\beta} (S^2), \qquad \beta \in (0,1),\\
X\supset D &=& C^{2+\beta}(S^2),
\end{eqnarray}
it turns out that 
\begin{eqnarray} \label{interpspaces}
D_A(\alpha, \infty)= C^{2\alpha+\beta} (S^2). 
\end{eqnarray}
This result is a straightforward adaptation of Theorem 3.1.12 in \cite{L1}, here with space variable $x\in S^2$ 
(or any compact manifold without boundary) instead of $x\in \mathbb{R}^n$.
\end{remark}
In order to state the strong stable manifold theorem, let us introduce the projection on the  unstable 
and strong stable eigenspaces respectively. More precisely: 
Let $\sigma(A)$ denote the spectrum of $A$ and take $\sigma_+(A)=\{\mu\in\sigma(A):\Real\mu\geq 0\}$ and 
 $\sigma_-(A)=\{\mu\in\sigma(A):\Real\mu < 0\}$. Note that $\sigma_+(A)$ consists of
finitely many eigenvalues of finite multiplicities associated with finitely many eigenfunctions.  In the case of the trivial equilibrium $v_*=0$, 
these are just the spherical harmonics.  
Take now $P_+$ to be the eigenprojection from $X$ onto the associated eigenspace and $P_-=(\id-P_+)$.
The following strong stable manifold theorem is an adaptation of \cite{L1}, 
Theorem 9.1.8 with the difference that we allow the linearization $A$ to have a zero eigenvalue. 
To compensate for this, we use the spectral gap: let $\omega_- = \sup \sigma_-(A)<0$, 
and $\eta \in (0,-\omega_-)$. For any time interval $I\subset \mathbb{R}$, 
we define the space of functions decaying at least as fast as $e^{-\eta t}$ by 
\begin{eqnarray*} \label{weightedspace}
C_{\eta}(I, D_A(\alpha+1, \infty))&=& \{ g: t\mapsto e^{\eta t} g(t) \in C(I, D_A(\alpha+1, \infty)) \},\\  
\Vert g \Vert_{C_{\eta}(I, D_A(\alpha+1, \infty)}&=&\sup_{t\in I} \Vert e^{\eta t}g \Vert_{C(I, D_A(\alpha+1, \infty)}.
\end{eqnarray*}
\begin{theorem}

For all $\alpha \in (0,1)$, there exists a constant $B_1>0 $ and a Lipschitz function $\psi$, 
\begin{equation} \label{strongstablepsi}
\psi: B(0,b_1) \subset P_-(C^{2\alpha+\beta}(S^2)) \rightarrow P_+(C^{2(\alpha+1)+\beta}(S^2)) , 
\end{equation}
differentiable at 0 with $\psi'(0)=0$, whose graph is the local strong-stable manifold for Equation~\eqref{eq:dyneq}. In other words:\\
For all $\tilde{v}_0 \in$ Graph$(\psi)$,  equation (\ref{eq:dyneq}) has  a unique solution 
$\tilde{v}\in C_{\eta}([0,\infty(,C^{2(\alpha+1)+\beta}(S^2))$ 
such that $\Vert \tilde{v}\Vert_{ C_{\eta}([0,\infty), C^{2(\alpha+1)+\beta}(S^2))} \leq B_1$.\\
And conversely there exists a constant $b_1>0$ such that if equation (\ref{eq:dyneq}) has a solution  
$\tilde{v}\in C_{\eta}([0,\infty), C^{2(\alpha+1)+\beta}(S^2))$ with 
$\Vert \tilde{v}\Vert_{ C_{\eta}([0,\infty),C^{2(\alpha+1)+\beta}(S^2))}\leq B_1$ 
and $\Vert P_-(\tilde{v}(0))\Vert_{ C^{2\alpha+\beta}(S^2)}\leq b_1$, 
then $\tilde{v}(0)\in$ Graph$(\psi)$. 
\end{theorem}

This is essentially part (ii) of Theorem 9.1.8 of~\cite{L1} or part (i) of Theorem 2.4 of \cite{L3}, i.e. the strongly stable 
manifold parts of those theorems.  Those theorems apply 
to a dynamical equation of the form~\eqref{eq:dyneq}, where $A$ is sectorial and the nonlinearity $Q$ may be regarded as a $C^{\infty}$ mapping from
a neighborhood 
$\mathcal{O}\subset D_{A}(\alpha+1,\infty)$ into $D_{A}(\alpha,\infty)$ and in addition satisfies $Q(0)=0, Q'(0)=0$.  Taking the interpolation 
spaces as generated by $C^{\beta}(S^2)$ and $C^{2+\beta}(S^2)$ as described above, these conditions hold in our case.  
However, the results from ~\cite{L1},
~\cite{L3} both assume that there is no null eigenvalue, which we do not assume.  It can nonetheless be checked that one retains a strongly 
stable manifold theorem even when an eigenvalue vanishes as long as one utilizes the spectral gap $(0,\omega_-)$ correctly.

In order to discuss the heteroclinics between the branches of equilibria (see Section 7, Remark \ref{heteroclinics}), 
we shall prove a center-unstable manifold theorem.  

To carry out this discussion we recall that 
the operator $A:D\to X$  is sectorial. 
Again, it is important that the nonlinearity $Q$ is a $C^{\infty}$ mapping $Q:C^{2(\alpha+1)+\beta}(S^2)\to C^{2\alpha+\beta}(S^2)$, 
which in addition satisfies $Q(0)=0, Q'(0)=0$.  Let $A_+=A_{P_+(C^{2\alpha+\beta}(S^2))}$ and 
$A_-=A_{P_-(C^{2(\alpha+1)+\beta}(S^2))}$.
With $x(t)=P_+\tilde{v}(t)$ and $y(t)=P_-\tilde{v}(t)$ we make the following decomposition  of our equation~(\ref{eq:v1.*}):       
\begin{align}
    x'(t)&=A_+x(t)+P_+Q(x(t)+y(t))\label{eq:system1}\\
    y'(t)&=A_-y(t)+P_-Q(x(t)+y(t))\label{eq:system2}
\end{align}

For the  center-unstable manifold theorem we shall first make a small modification of the 
system~(\ref{eq:system1})-(\ref{eq:system2}).  Let $\rho:P_+(C^{2\alpha+\beta}(S^2)) \to \R$ be a $C^{\infty}$ cutoff
function such that 
\[
    0\leq \rho(x)\leq 1,\,\, \rho(x)=1 \,\, \text{if}\,\,  \norm{x}_{C^{2\alpha+\beta}(S^2)}\leq 1/2,\,\, \rho(x)=0 \,\, \text{if}\,\, 
    \norm{x}_{C^{2\alpha+\beta}(S^2)}\geq 1.
\]
The existence of such a function is no problem since $P_+(C^{2\alpha+\beta}(S^2))$ is finite dimensional.  
For $q>0$ small we consider instead the  system modified by cutoff
\begin{align}
       x'(t)&=A_+x(t)+Q_+(x(t),y(t))\label{eq:system1'}\\
    y'(t)&=A_-y(t)+Q_-(x(t),y(t)) \label{eq:system2'},
\end{align}
where  $Q_+(x,y)=P_+Q(\rho(x/q)+y)$ and $Q_-(x,y)=P_-Q(\rho(x/q)x+y)$, and $q$ is called the cut-off parameter.

note that the above construction can be carried out in such a way 
that it is equivariant with respect to $O(3)$ in the case that the original equilibrium is trivial $v_*\equiv 1$.  
Indeed, it is certainly the case that our original equation, as well as the spaces
$X_+=P_+(C^{2,\beta}(S^2))$, $X_-=P_-(C^{2,\beta}(S^2))$,
are equivariant under this group action.  Hence we must only check that we can arrange for the modified 
equations~(\ref{eq:system1'})-(\ref{eq:system2'}) to be equivariant as well.  
In order to do this, we use an isometry (using the $L^2$ norm on $X_+$) $F:X_+:\to\R^m$ 
to identify $X_+$ with $\R^m$.
Thus, arranging that the cutoff function $\rho$ is such that $\rho\circ F$ is spherically symmetric, we can ensure that the 
modified equations~(\ref{eq:system1'})-(\ref{eq:system2'}) are $O(3)$ equivariant.  
\begin{theorem}
Let now $A$ and $Q$ be as in the first paragraph of this section and consider equation~(\ref{eq:dyneq}) 
as well as  the modified equations ~(\ref{eq:system1'})-(\ref{eq:system2'})  constructed with a cutoff 
$\rho$ as in the remarks of the previous paragraph. \\
In particular $A$ is sectorial and admits finitely many unstable eigenvalues of finite multiplicities. 
The nonlinearity $Q$ is assumed to be a mapping  $C^{\infty}(\mathcal{O}, C^{2\alpha+\beta}(S^2))$, where 
$\mathcal{O}$ is a neighborhood of the origin in $C^{2(\alpha+1)+\beta}(S^2)$. \\
Then there exists a cutoff parameter $q_1$ such that if $q<q_1$, then there is a continuously differentiable function 
$\Xi:P_+(C^{2\alpha+\beta}(S^2)) \to P_-C^{2(\alpha+1)+\beta}(S^2))$ 
whose graph is invariant for the system~(\ref{eq:system1'})-(\ref{eq:system2'}).  Graph$\Xi$ is called a local center-unstable manifold.
In addition
\[
      \Xi'(x)(A_+x+Q_+(x,\Xi(x)))=A_-\Xi(x)+Q_-(x,\Xi(x)), x\in P_+(X).
\]
In the case of the trivial equilibrium,
$\Xi$ is equivariant with respect to $O(3)$, and so is its graph.\\
The function $\Xi$ can be made $C^k$ by choosing a possibly smaller bound $q_k$ for the cutoff parameter. 
\end{theorem}
This theorem can be found in~\cite{L1}.  For a similar theorem see the nice succinct work of Mielke~\cite{Mielke}.

Of course the finite dimensional center-unstable manifold contains a center manifold which 
is obtained by using the projectors $P_0$ and $P_h$ onto the center eigenspace associated to the 
zero eigenvalue and the ``hyperbolic eigenspace" of infinite dimension associated to the nonzero eigenvalues respectively. 
With $x(t)=P_0\tilde{v}(t)$ and $y(t)=P_h\tilde{v}(t)$ we make the following decomposition  of our equation~(\ref{eq:v1.*}):       
\begin{align}
    x'(t)&=P_0Q(x(t)+y(t))\label{eq:systemcenter1}\\
    y'(t)&=A_hy(t)+P_hQ(x(t)+y(t))\label{eq:systemcenter2}
\end{align}
This proves the following theorem. 
\begin{theorem}
\label{centermnf}
Let $\mathcal{O}$ be a neighborhood of the origin in $C^{2(\alpha+1)+\beta}(S^2)$ 
in which the nonlinearity $Q$ is continuously differentiable, $Q\in C^1(\mathcal{O}, C^{2(\alpha+1)+\beta}(S^2) )$. 
The linear part $A$ is sectorial and admits only a finite number of nonnegative eigenvalues of finite multiplicity, $Q(0)=Q'(0)=0$. \\
Then there exists a cut-off parameter $q_1$ such that for $0<q\leq q_1$, there exists a Lipschitz continuous function 
$\Phi: P_0(C^{2\alpha}(S^2))\rightarrow P_h (C^{2(\alpha+1)}(S^2))$ whose Graph is invariant under the flow of the system modified by cut-off
\begin{align}
       x'(t)&=Q_0(x(t),y(t))\label{eq:cutoffcenter1}\\
    y'(t)&=A_hy(t)+Q_h(x(t),y(t)) \label{eq:cutoffcenter2},
\end{align}
where  $Q_0(x,y)=P_0Q(\rho(x/q)+y)$ and $Q_h(x,y)=P_hQ(\rho(x/q)x+y)$.

Graph$\Phi$ is called a (local) center manifold. 
\end{theorem}

\section{Anisotropic self-similar blow-up}
\label{blowup}
\setcounter{equation}{0}

In this section we apply the equivariant branching Lemma \ref{lm3.1}
to obtain anisotropic self-similar blow-up solutions
\begin{equation}
 \label{eq5.1} 
u(r,p) =\left( \frac{\lambda}{2} \left( \frac{1}{r}-1\right) \right)^{-\frac{1}{2}} (v+1)
\end{equation}
with isotropically prescribed scalar curvature $R(r)=(\lambda +2)/r^2$, as we announced in the
introduction. By (\ref{eq:v1.*}) this requires $v=v(p)$ to be a nonhomogeneous
solution of the equilibrium equation
\begin{equation}
0=\Delta v + \lambda f(v)\label{eq5.2}
\end{equation}
for the blow-up shape $v(p), \ p\in S^2$. We also recall the
nonlinearity 
\begin{equation}
f(v)  v - \frac12 v^2 /(v+1) \label{eq5.3}
\end{equation}
and the coefficient
\begin{equation}
\lambda = r^2 R(r)-2
 \label{eq5.4}
\end{equation}
Anisotropic, i.e. nonhomogeneous, solutions $v$ of (\ref{eq5.2}) will
arise by $O(3)$ symmetry breaking bifurcations with maximal
isotropy and one-dimensional fix spaces, as indicated in Lemma
\ref{lm2.1}. 

To apply the equivariant branching Lemma \ref{lm3.1} we first rewrite
(\ref{eq5.2}) as 
\begin{equation}
0= G(\lambda,v):= v + (\Delta-1)^{-1}(v+\lambda f(v)). \label{eq5.5}
\end{equation}
For Banach spaces $X=Y$ we could choose the Hilbert space setting $X=
L^2(S^2)$. For technical consistency with our dynamics approach of
H\"older type, above, we work in
\begin{equation}
X=Y:= C^{\beta}(S^2)\label{eq5.6}
\end{equation}
instead.

We notice immediately, that $G$ is a compact perturbation of
identity. Indeed the solution $(\Delta-1)^{-1}g$ of the linear Poisson
equation $\Delta v -v = g$ on $S^2$ is a bounded linear operator
\begin{equation}
(\Delta-1)^{-1}: \ C^{\beta}(S^2) \to C^{\beta}(S^2),\label{eq5.7}
\end{equation}
and the embedding $C^{2+\beta}(S^2) \hookrightarrow
C^{\beta}(S^2)$ is compact by Arzela-Ascoli.

Obviously $G$ is $C^k$, for $1+v >0$ strictly positive and any $k \geq
0$, by analyticity of $f$ on this subdomain of $\mathbb{R}$, and of
$C^{0,\alpha}$. A trivial branch $G(\lambda,0) = 0$ of homogeneous,
$O(3)$-fixed solutions arises from $f(0) = 0$. Orthogonal invariance
of the Dirichlet integral $\frac12 \int|\nabla v|^2 dp$
under the representation $(\varrho(\gamma)v)(p):= v(\gamma^{-1}p)$ of
(\ref{eq2.8}) for $\gamma \in \Gamma:= O(3)$ implies
$O(3)$-equivariance (\ref{eq2.3}) of $(\Delta-1)^{-1}$ on $ v\in
C^{0,\alpha}(S^2)$. Obviously $\varrho(\gamma)$ also commutes with the
point evaluation $v\mapsto f(v)$ in $C^{\beta}(S^2)$. This proves
$\Gamma$-equivariance of $G$ and establishes assumptions
(\ref{eq3.1})--(\ref{eq3.3}) of the equivariant branching Lemma
\ref{lm3.1}.

Consider the $O(3)$-equivariant linearization
\begin{eqnarray}
L(\lambda)v&:=& D_v G(\lambda,0)v = v + (\Delta-1)^{-1}(v+\lambda
f'(0)v) =\label{eq5.8}\\[1mm]
&\ =& v+(1+\lambda) (\Delta-1)^{-1}v\nonumber
\end{eqnarray}
at the trivial branch $G(\lambda,0) = 0$ next. Again $L(\lambda)$
is a compact linear perturbation of identity on $X=Y=
C^{\beta}(S^2)$. Hence $L(\lambda)$ is Fredholm of index zero on
$X$, and on any closed linear and $L(\lambda)$-invariant subspace of
$X$. In particular $L^K(\lambda)$, the restriction of $L(\lambda)$ to
the subspace $X^K ={\rm Fix}(K)$ of $K$-fixed elements in $X$, is Fredholm
for any subgroup $K$ of $\Gamma=O(3)$.

By standard spectral theory of the Laplace-Beltrami operator $\Delta$
on the 2-sphere $S^2$ we obtain nontrivial kernel
\begin{equation}
\{0\} \neq \ker L(\lambda) = V_\ell\label{eq5.9}
\end{equation}
if and only if, 
\begin{equation}
\lambda = \lambda_\ell := \ell(\ell+1)\label{eq5.10}
\end{equation}
for some $\ell = 0,1,\ldots$. The resulting kernels are the
$L^2$-complete spherical harmonics subspaces $V_\ell$ of
$X=C^{0,\alpha}(S^2)$; see (\ref{eq2.9}), (\ref{eq2.10a}) in Section
\ref{sym}.

To check the transversality condition (\ref{eq3.6}) of the equivariant
branching Lemma \ref{lm3.1} we calculate
\begin{equation}
D_\lambda L(\lambda) \mathbf{e} = (\Delta-1)^{-1} \mathbf{e} = \frac{1}{1+\ell(\ell+1)}
\mathbf{e}. \label{eq5.11}
\end{equation}
By $L^2$ self-adjointness of $\Delta$, however, $\mathbf{e} \in \ker L(\lambda)$
is $L^2$ orthogonal to   range $L(\lambda)$ and, a fortiori, orthogonal to
any subspace  range $L^K(\lambda)$. This proves transversality
condition (\ref{eq3.6}).

\begin{theorem}\label{th5.1}
{~}

Let $\ell=1,2,3,\ldots$ and let $K < O(3)$ be an isotropy subgroup of
$O(3)$ such that 
\begin{equation}
V^K_\ell = {\rm Fix}(K) = \ {\rm span}\ \{\mathbf{e}\} \label{eq5.12}
\end{equation}
is one-dimensional, for the representation $(\varrho(\gamma)v)(p):=
v(\gamma^{-1}p)$ of $O(3)$ in the spherical harmonics space $V_\ell$ of
dimension $2 \ell+1$; see (\ref{eq2.8}), (\ref{eq2.10a}).

Then there exists a unique smooth local curve
\begin{equation}
s \mapsto (\lambda(s),v(s)) \in \mathbb{R} \times C^{\beta}(S^2)\label{eq5.13}
\end{equation}
of solutions $(\lambda,v)$ of the self-similar blow-up equation
(\ref{eq5.2}) such that
\begin{equation}
\lambda(0)=\lambda_\ell = \ell(\ell+1), \ v(0) = 0, \ v'(0) = \mathbf{e} \label{eq5.14}
\end{equation}
and the isotropy  of $v(s)$ is exactly $K < O(3)$, for $s \neq 0$:
\begin{equation}
\Gamma_{v(s)} = K \neq O(3).\label{eq5.15}
\end{equation}
In  particular $v(s)$ is nonhomogeneous and hence anisotropic, for
$s\neq 0$.
\end{theorem}

{\bf Proof:}
{~}

By the arguments above and assumption (\ref{eq5.15}), the equivariant
branching Lemma \ref{lm3.1} asserts the existence of a local branch
$(\lambda(s),v(s)) \in \mathbb{R} \times X^K$, of class $C^k$ for any
$ k \geq 1$. Smoothness follows, for the $k=1$ branch and any $s\neq
0$, by direct implicit function theorem. At $s=0$ smoothness follows
by uniqueness of the branches, even if the $s$-neighborhoods had shrunk
to zero for $k \to +\infty$.

To check for the precise isotropy $\Gamma_{v(s)} = K$ is a little more
subtle. In Section \ref{sym} we have already remarked that isotropies $K$ with
one-dimensional fix spaces are maximal in their irreducible
representation. This proves
\begin{equation}
\Gamma_{\mathbf{e}} = K \neq O(3)\label{eq5.16}
\end{equation}
for the eigenvector $\mathbf{e}$ which spans $V^K_\ell$.

To extend this observation from $\mathbf{e}$ to 
\begin{equation}
v(s) = s\mathbf{e} + o(s),\label{eq5.17}
\end{equation}
for small enough $s$, we argue indirectly: suppose there exists a
sequence $0\neq s_i \to 0$ such that

\begin{equation}
K < \Gamma_{v(s_i)} \le \Gamma = O(3). \label{eq5.18}
\end{equation}
In $\Gamma=O(3)$ this leaves only finitely many choices for
$\Gamma_{v(s_{i})}$. Indeed the only infinite ascending series of
closed subgroups are generated by the rotations $\mathbb{Z}_k$ over
$2\pi/k$ in $SO(2) \le SO(3)$. This case is easily excluded, by
direct inspection of Fourier series in the azimuthal latitude
$\varphi$. Without loss of generality we may therefore assume
\begin{equation}
K < \Gamma_{v(s_{i})} \equiv \Gamma_\infty \le \Gamma = O(3)\label{eq5.19} 
\end{equation} 
does not depend on $i$.

Since $v'(0) =\mathbf{e}$ we have $v(s_{i})/s_i \to \mathbf{e}$ for $0 \neq s_i \to
0$. Continuity of the representation $\varrho$ therefore implies

\begin{equation}
\Gamma_{v(s_{i})} = \Gamma_\infty \le \Gamma_{\mathbf{e}} = K
\end{equation}
which contradicts our indirect assumption (\ref{eq5.18}). This
argument, which is a special variant of the virtual isotropy
proposition in \cite{FiedlerHabil}, completes the proof of (\ref{eq5.15})
and of the theorem.

\hfill$\bowtie$

We now apply Theorem \ref{th5.1} on $O(3)$ symmetry breaking of the
trivial equilibrium $v=0$ in (\ref{eq5.5}) to the maximal isotropy subgroups
$K$ of the spherical harmonics representations of Lemma
\ref{lm2.1}. We address the {\it axisymmetric cases} $K \geq SO(2)$,
for all $\ell$, and the {\it discrete isotropies} $K$, for $\ell =3,4$,
separately, in Corollaries \ref{cor5.2} and \ref{cor5.3} below.

\begin{corollary}\label{cor5.2}
{~}

Let $\ell > 0$ and consider axisymmetric fix-groups
\begin{equation}\label{eq5.21}
K: = \left\{
\begin{array}{ll} 
O(2)^- & {\rm for \ odd} \ \ell,\\[2mm]
O(2) \oplus \mathbb{Z}^c_2 & {\rm for \ even} \ \ell.
\end{array}\right.
\end{equation}

Then the conclusions of Theorem \ref{th5.1} hold and provide a unique
smooth local branch $(\lambda(s),v(s))$ of nontrivial solutions
bifurcating from $(\lambda(0),v(0)) = (\lambda_\ell,0)$ at
$\lambda_\ell = \ell(\ell+1)$. The solutions  $v(s)$ at $\lambda=\lambda(s)$ possess
isotropy $\Gamma_{v(s)} = K$ for $s\neq 0$. In particular they are
nonhomogeneous but {\it axisymmetric}, i.e. independent of the
azimuthal angle $\varphi$ on $S^2$.
\end{corollary}

{\bf Proof:}
{~}

Lemma \ref{lm2.1} implies that the one remaining assumption
(\ref{eq5.12}) of Theorem \ref{th5.1} holds for $K$ from
(\ref{eq5.21}):
\begin{equation}
 \dim {\rm Fix}\ (K) = \dim V^K_\ell = 1. \label{eq5.22}
\end{equation}
This proves the corollary.

\hfill $\bowtie$

\begin{corollary}\label{cor5.3}
{~}

Let $\ell=3$ and $K\in \{O^-,D^d_6\}$, or let $\ell=4$ and $K= O
\oplus \mathbb{Z}^c_2$.

Then the conclusions of Theorem \ref{th5.1} hold and provide a unique
smooth local branch $(\lambda(s),v(s))$ of nontrivial solutions
bifurcating from $(\lambda(0),v(0)) = (\lambda_\ell,0)$ at
$\lambda_\ell = \ell(\ell+1)$. The solutions for $v(s)$ possess
isotropy $\Gamma_{v(s)} = K$, for $s \neq 0$. In particular they are
nonhomogeneous with hexagonal or octahedral   symmetry as given by $K$. 
\end{corollary}

{\bf Proof:}
{~}

The proof of the previous corollary applies verbatim.

\hfill$\bowtie$

In conclusion we obtain self-similar but anisotropic blow-up in the isotropic
scalar curvature equation, with remaining octahedral, hexagonal, or
axisymmetry as specified by the isotropy $K$ above.

Of course our result applies equally well to all maximal isotropy
subgroups $K$ of $O(3)$ on $V_\ell$, as long as the spaces Fix$(K)$ of
$K$-fixed vectors in $V_\ell$ remain real one-dimensional. See \cite{CLM, LautHabil} 
for lists of these cases.

%
%
%

\section{Anisotropic asymptotically self-similar blow-up}
\label{asym}
\setcounter{equation}{0}

In this section we address dynamic aspects of the self-similarly
rescaled scalar curvature blow-up problem
\begin{equation}
(1+v)^{-2}\partial_t v = \Delta v + \lambda f(v);\label{eq7.1}
\end{equation}
see (\ref{eq:v1.*}). More specifically we study the Morse indices of the $O(3)$
symmetry breaking bifurcations $(\lambda(s),v(s))$ from the trivial
homogeneous solution branch $(\lambda,v)=(\lambda,0)$ at the
bifurcation points $(\lambda(0),c(0)) = (\lambda_\ell,0)$ with
\begin{equation}
\lambda_\ell = \ell(\ell-1)\label{eq7.2}
\end{equation}
for $\ell=1,2,3,\ldots $. The nonlinearity $f(v)$ takes the specific
form
\begin{equation}
f(v) = v- \frac12  v^2/(1+v)
 = v - \frac12 v^2 + \frac12
v^3\mp \ldots 
\label{eq7.3}
\end{equation}
In particular the first Taylor coefficients are 
\begin{equation}
f(0) = 0, \ f'(0) = 1, \ f''(0) = -1, \ f'''(0)= 3.\label{eq7.4}
\end{equation}
We give a complete account of the low-dimensional bifurcation diagrams
$\ell=1,2,3,4$ in Theorem \ref{th7.3} based on Corollary \ref{cor5.2}
and \ref{cor5.3}. We omit the exceptionally trivial case $\ell=0, \
\lambda_\ell =0$, which is linear and only features vertical
bifurcation of spatially homogeneous, isotropic equilibria $v \equiv
\ {\rm const}$.

Our local analysis starts from the standard equivariant center manifold
theorem 
with
equivariance group $\Gamma = O(3)$ for the Morse equivalent but semilinear parabolic PDE 
\begin{equation}
\label{6.4a}
\D_{t}v= \Delta v+\lambda f(v).
\end{equation}
Adding the equation $\dot{\lambda}
= 0$, artificially, we obtain a trivially extended semiflow
$S^\lambda(t)$ of class $C^k$ on the extended phase space
$\mathbb{R} \times \left( X \cap \{1+v > 0\} \right)$.  At bifurcation
points $(\lambda,0) = (\lambda_\ell,0), \ \lambda_\ell =
\ell(\ell+1)$, we immediately obtain a reduced flow of (\ref{6.4a}) on the center manifold tangent to $V_{\ell}$ in
parametrized and equivariant form:
\begin{equation}
\dot{v}_c = \Phi(\lambda,v_c),\qquad v_c \in V_\ell. \label{eq7.5}
\end{equation}

In Lemma \ref{lm7.1} we first compare this flow, truncated to order 2
or 3, with the detailed discussion of the generic $O(3)$-equivariant
bifurcations on $V_\ell$ in \cite{LautHabil, CLM}. Specifically, we establish
nondegenerate transcriticality and the sub- or supercritical
nondegenerate pitchfork type of the maximal isotropy branches of
Corollaries \ref{cor5.2}, \ref{cor5.3}, respectively. We call the
bifurcation of $(\lambda(s),v(s))$ at $s=0, \ (\lambda(0),v(0)) =
(\lambda_\ell,0)$ nondegenerate and
\begin{equation}
\label{eq7.6}
\begin{array}{lll}
{\it transcritical}, & {\rm if} & \lambda'(0) \neq 0\\[2mm]
{\it pitchfork}, & {\rm if} & \lambda'(0) = 0 \neq \lambda''(0).
\end{array}
\end{equation}
Because these branches of equilibria must reside in the local center
manifold of (\ref{6.4a}), they establish nondegeneracy of the respective quadratic and
cubic equivariants of $x_c \mapsto \Phi(\lambda_\ell,x_c)$, and thus
determine the resulting unstable dimensions of the bifurcating
equilibria $v(s)$ inside the center manifold. This little trick avoids very messy direct computations of, and inside, the center manifolds. 

In a second step, in Proposition \ref{prop7.2}, we then establish normal
hyperbolicity of the associated group orbits $\Gamma v(s) \cong
\Gamma/\Gamma_{v(s)}$ of equilibria $v_*$. We determine the codimensions of
the local strong stable manifolds $W^{ss} ={\rm graph} \Psi_{ss}$, for each of the
symmetry breaking branches specified in Theorems \ref{th5.1}, \ref{th7.3} and
Corollaries \ref{cor5.2}, \ref{cor5.3}. Since 
\begin{equation}
v(\cdot,t) \to v_* \label{eq7.7}
\end{equation}
converges exponentially in $C^{\beta}$, on graph $\Psi_{ss}$, the
solutions in graph $\Psi_{ss}$ provide anisotropic and only
asymptotically self-similar solutions $u$ with 
blow-up at $r=1$ of the isotropic scalar curvature $R(r)=(\lambda+2)/r^2$.
\begin{lemma}\label{lm7.1}
{~}

Consider the symmetry breaking branches $(\lambda_K(s),v_K(s))$ of
equilibria of Equation {\rm (\ref{eq7.1}), \ref{6.4a}} with isotropy $K$ which bifurcate
from $(\lambda_K(0),v_K(0))=(\lambda_\ell,0), \ \lambda_\ell =
\ell(\ell+1), \ \ell=1,2,3,\ldots,$ as established in Corollaries
$\ref{cor5.2}$ and $\ref{cor5.3}$.

Then the bifurcations of the axisymmetric solutions $K \geq SO(2)$ of
Corollary $\ref{cor5.2}$ are nondegenerate transcritical, for all even
$\ell > 0$. The remaining solution branches at $\ell = 3,4$
established in Corollary $\ref{cor5.3}$ are 
\begin{equation}\label{branches}
\begin{array}{lll}
{\rm nondegenerate} & {\rm transcritical},& {\rm for} \  \ell \ {\rm even}\\[2mm]
{\rm nondegenerate} & {\rm pitchfork}, &{\rm for}  \ \ell \ {\rm odd}
\end{array}
\end{equation}

At $\ell = 3,4$ the pitchfork curvatures $\lambda''_K(0)$ of the
branches with maximal isotropy $K$ are given in Table \ref{tableofnonvanishing}.
Here we use $v'_K(0) = \mathbf{e}_K$ to normalize the parametrization of all
branches $(\lambda_K(s),v_K(s))$ by $s$ such that

\begin{equation}
\langle \mathbf{e}_K, v_K(s) \rangle = s\label{eq7.8}
\end{equation}
holds for the $L^2$ scalar product $\langle\cdot,\cdot\rangle$ on
$S^2$ with the normalized eigenvector $\mathbf{e}_K$ of the space $V^K_\ell$ of $K$-fixed
vectors in $V_\ell$. 
\end{lemma}
\begin{table}
\begin{center}
\begin{tabular}{||c|c||c|c|c|c|c||} \hline\hline
 \multicolumn{2}{||c||}{}&&&&&\\
 \multicolumn{2}{||c||}{isotropy $K$} & $O(2)^-$ & $O(2) \oplus \mathbb{Z}^c_2$ &$O^-$ & $O
\oplus \mathbb{Z}^c_2$  &$D^d_6$ \\ 
 \multicolumn{2}{||c||}{}&&&&&\\
\hline
&&&&&&\\
&$\mathbf{e}_K$ &$Y_{10}$& - & - & - & - \\
&&&&&&\\
 \cline{2-7} 
 $\ell=1$ &&&&&&\\
& $\lambda''(0)$&$-\frac{3}{5 \pi}$& - & - & - & - \\
&&&&&&\\
\hline
&&&&&&\\
 &$\mathbf{e}_K$& - &$Y_{20}$& - & - & - \\
 &&&&&&\\
 \cline{2-7} 
 $\ell=2$&&&&&&\\
& $\lambda'(0)$ & - &$\frac{\sqrt{5}}{14 \sqrt{ \pi}}$&  - & - & -\\
&&&&&&\\
\hline
&&&&&&\\
 &$\mathbf{e}_K$ &$Y_{30}$& - &$\sqrt{2}$Re$Y_{32}$& - & $\sqrt{2}$Re$Y_{33}$\\
&&&&&&\\
 \cline{2-7} 
 $\ell=3$&&&&&&\\
& $\lambda''(0)$&$-\frac{2954}{715 \pi}$& - &$-\frac{1050}{143 \pi}$& - &$-\frac{665}{286 \pi}$\\
&&&&&&\\
\hline
&&&&&&\\
& $\mathbf{e}_K$& -  &$Y_{40}$& - &$\frac12 \sqrt{\frac{7}{3}} \left( Y_{40}+ 2\sqrt{\frac{5}{14}} \rm{Re}Y_{44} \right) $& - \\
&&&&&&\\
 \cline{2-7} 
 $\ell=4$&&&&&&\\
& $\lambda'(0)$& - &$\frac{243}{2002 \sqrt{\pi}}$& - &$\frac{9\sqrt{21}}{286 \sqrt{\pi}}$& - \\
&&&&&&\\
\hline\hline
\end{tabular}
\end{center}\vspace{0.3cm}
\caption{Exact values of normalized eigenvectors $\mathbf{e}_K$ spanning isotropy fix spaces $V_{\ell}^K$
and of first nonvanishing derivatives $\lambda'_K(0), \lambda''_K(0)$ for
symmetry breaking equilibrium branches $(\lambda_K(s),v_K(s))$ of
(\ref{eq7.1}), (\ref{6.4a}), with maximal isotropy $K$. \label{tableofnonvanishing}}
\end{table}


{\bf Proof:}
{~}

We first calculate the derivatives of $\lambda(s)$ at $s=0$, in the
general setting of the symmetry breaking bifurcation Theorem
\ref{th5.1}. We then specify to the axisymmetric cases $K \geq
SO(2)$, and to the maximal isotropies $K$ in the $V_\ell$
representations of $\Gamma = O(3)$.

We recall the general setting
\begin{eqnarray}
G(\lambda,v) &=& v + (\Delta-1)^{-1}(v+\lambda f(v)) =
0\nonumber\\[1mm]
L(\lambda)v= \partial_v G(\lambda,0)v &=& v +
(\Delta-1)^{-1}(1+\lambda)v\label{eq7.9}\\
f(v) &=& v-\frac12 v^2(v+1)^{-1} = v - \frac12 v^2 + \frac12
v^3\mp \ldots\nonumber
\end{eqnarray}
from (\ref{eq5.3}), (\ref{eq5.5}), (\ref{eq5.8}). By definition 
\begin{equation}
{\rm Fix}\ (K) \cap \ker L(\lambda_\ell) = \ {\rm Fix}\ (K)\cap V_\ell =
\ {\rm span} \ \{\mathbf{e}_K\}, \label{eq7.10}
\end{equation}
where we choose $\mathbf{e}_K$ to be an $L^2$-normalized unit vector. See
Theorem \ref{th5.1} and(\ref{eq5.12})--(\ref{eq5.14}), where
dependence on the chosen isotropy $K$ was suppressed. For the explicit
entries $\mathbf{e}_K$ in Table \ref{tableofnonvanishing} see \cite{CLM}  and \cite{FiedMisch}. The case
$\mathbf{e}_K = Y_{\ell 0}$ for $K \geq SO(2)$ is obvious: the fix space
generator $\mathbf{e}_K \in V_\ell$ is independent of the azimuthal angle
$\varphi$ and hence normalized to $\pm Y_{\ell 0}(\vartheta)$ in the
spherical harmonics space $V_\ell$.

From Section \ref{sym} we recall that $- {\rm id} \in \mathbb{Z}^c_2$ acts on
$V_\ell$ by multiplication with $(-1)^\ell$. In particular $- {\rm id}
\notin K$ for any nontrivial isotropy $K$ on $V_\ell$ and odd $\ell$.
Conjugation by $- {\rm id} \in \mathbb{Z}^c_2$, on  the other hand,
leaves $K$ invariant. Uniqueness of the bifurcation branch
$(\lambda_K(s),v(s))$ and the normalization (\ref{eq7.8}) therefore
imply

\begin{equation}
\lambda_K(-s) = \lambda_K(s), \quad v_K(-s) = -v_K(s)\label{eq7.11}
\end{equation}
for any isotropy $K$ and odd $\ell$. In particular
\begin{equation}
\lambda'_K(0) = 0\label{eq7.12}
\end{equation} 
and all such symmetry breaking bifurcations are pitchfork candidates.

To calculate the first derivative of $\lambda(s)$ at $s=0$, in
general, we suppress $K$ again and differentiate $G(\lambda(s),v(s)) = 0$
twice with respect to $s$:
\begin{eqnarray}
0 &=& \partial_\lambda G\lambda' + \partial_v G v'\label{eq7.14}\\[3mm]
0 &=& \partial_\lambda G\lambda'' + \partial^2_\lambda G(\lambda')^2 + 2 \partial_\lambda
\partial_v G \lambda' v' + \partial^2_v G(v')^2 + \partial_v Gv''.\label{eq7.15}
\end{eqnarray}
At $s=0$, $v(0) = 0$ all $\lambda$-derivatives of $G$ vanish because
$G(\lambda,0) \equiv 0$ along the trivial branch. Self-adjointness of
$\Delta$ implies orthogonality 
\begin{equation}
\langle \mathbf{e}, \partial_v G\rangle = 0\label{eq7.13}
\end{equation}
in $L^2(S^2)$ at $s=0$ where $\partial_v G= L(\lambda_\ell)$. Taking the
scalar product of (\ref{eq7.12}) with $\mathbf{e}$ at $s=0$ implies
\begin{equation}
0=2 \lambda'(0) \langle \mathbf{e},(\Delta-1)^{-1}f'(0)\mathbf{e}\rangle + \langle
\mathbf{e},(\Delta-1)^{-1} f''(0)\mathbf{e}^2\rangle. \label{formula}
\end{equation}
With $\Delta \mathbf{e} = - \lambda_\ell \mathbf{e}, \ f'(0) = 1, \ f''(0)=-1$ and the
$L^2$-normalization $\langle \mathbf{e},\mathbf{e} \rangle = 1$ we obtain
\begin{equation}
\lambda'(0) = \frac12 \langle \mathbf{e},\mathbf{e}^2\rangle . \label{lambda'0}
\end{equation}

We claim next that 
\begin{equation}
\lambda'_K(0) \neq 0 \label{eq7.16}
\end{equation}
is nondegenerate transcritical, for all axisymmetric solutions $K =
O(2) \oplus  \mathbb{Z}^c_2$ and even $\ell = 2,4,6,\ldots$. We use \cite{AbraSte}, 27.9.1
 to explicitly calculate $\langle
\mathbf{e},\mathbf{e}^2\rangle = \int\limits_{S^2}(Y_{\ell 0})^3$ as follows
\begin{eqnarray}
\langle \mathbf{e},\mathbf{e}^2\rangle &=& \sqrt{\frac{(2\ell+1)^3}{4\pi}}  \binom{\ell\ \ell\ \ell} {0\ 0\ 0}^2=\nonumber\\
&=& \sqrt{\frac{2\ell+1}{4\pi}} (\ell\ \ell\ 0\ 0|\ell\ \ell\ \ell\
 0)^2 =\nonumber\\[-2mm]
&&\label{eq7.17}\\ [-2mm]
&=& \sqrt{\frac{(2\ell+1)^3}{4\pi}} \ \frac{(\ell !)^3}{(3\ell+1)!}
 \left(\sum^\ell_{k=0}(-1)^k \left(\frac{\ell !}{k!(\ell-k)!}\right)^3\right)^2 =\nonumber\\
&=& \sqrt{\frac{(2\ell+1)^3}{4\pi}} \ \frac{(\ell !)^3}{(3\ell+1)!}
 \left(2 \sum\limits^{\ell/2}_{k=0} (-1)^k \binom{\ell}{k}^3\right)^2 \neq
 0. \nonumber
\end{eqnarray}
We have used $\ell$ even in the last equality; for odd $\ell$ we
encounter a rather circuitous confirmation of $\lambda'(0) = 0$. The
alternating sum from $k=0$ to $k=\ell/2$ is nonzero because the binomial
coefficients $\binom{\ell}{k}$ increase strictly for
$k=0,\ldots,\ell/2$. This proves nondegenerate transcriticality
(\ref{eq7.16}) for axisymmetric solutions and even $\ell$. For example
we obtain the entry for $\ell=4, K=O(2) \oplus \mathbb{Z}^c_2$ in Table
\ref{tableofnonvanishing}.

We consider the octahedral case, $\ell=4$, $K= O \oplus \mathbb{Z}^c_2, \ \mathbf{e}=\mathbf{e}_K =
\frac 12 \sqrt{7/3}\ \bigg(Y_{40} + 2 \sqrt{5/14} \ {\rm Re}\
Y_{44}\bigg)$ next. Up to positive coefficients we obtain
\begin{equation}
\langle \mathbf{e},\mathbf{e}^2\rangle = \int\limits_{S^2} \bigg(\frac12\sqrt{\frac73}
\bigg(Y_{40} +
\sqrt{\frac{5}{14}}\bigg(Y_{44}+Y_{4,-4}\bigg)\bigg)\bigg)^3 =
\frac{9}{143} \sqrt{\frac{21}{\pi}} \neq 0. \label{eq7.18} 
\end{equation}
The last equality can be verified, e.g., by Mathematica or similar
programs. Alternatively it follows by hand, again, using \cite{AbraSte}, 27.9.1. 

It remains to determine the three pitchfork curvatures
$\lambda''_K(0)$ for $\ell=3$,  $K\in\{O(2)^-,O^-,D^d_6\}$, $\mathbf{e} =
\mathbf{e}_K$. A similarly tedious but crucial calculation was already
performed in \cite{FiedMisch}, based on the generators $\mathbf{e}=\mathbf{e}_K$ of
the isotropy fix spaces $V^K_3$ from \cite{LautHabil, CLM}. We
differentiate the second derivative (\ref{eq7.15}), once again, with
respect to $s$ at $s=0$, using $\lambda'(0) = 0$ and $\partial^k_\lambda$ G=0
for all $k$:
\begin{equation}
0=3 \partial_\lambda \partial_v G \lambda'' v' + \partial^3_v G(v')^3 + 3\partial^2_v G v' v'' +
\partial_v Gv'''. \label{eq7.19}
\end{equation}

We test against the $L^2$ unit vector $\mathbf{e}=v'(0)$ as before, and insert
the derivatives of $G,f$ to obtain
\begin{equation}
0= 3 \lambda'' + 3 \lambda_\ell \langle \mathbf{e},\mathbf{e}^3\rangle - 3 \lambda_\ell
\langle \mathbf{e}^2,v''\rangle \label{eq7.20}
\end{equation}
at $s=0$. Evaluating (\ref{eq7.15}), once again with $\partial_{\lambda}G=0$, $\lambda'=0$, $v'=\mathbf{e}$, at $s=0$ we observe
\begin{equation}
(\Delta + \lambda_\ell) v''(0) = + \lambda_\ell \mathbf{e}^2. \label{eq7.21}
\end{equation}
Because $v''(0) \in X^K$ is also $L^2$-orthogonal to the
one-dimensional kernel $V^K_\ell =$ span $\{\mathbf{e}\}$ of $(\Delta +
\lambda_\ell)$ in $X^K$, by normalization (\ref{eq7.8}), we can write
$v''(0)$ uniquely as a finite sum
\begin{equation}
v''(0) = \sum\limits_{\underset{|m'| \le \ell'}{\ell'\neq\ell}} \ 
\frac{\lambda_\ell}{\lambda_\ell-\lambda_{\ell'}}\ c_{\ell'm'}\ Y_{\ell'm'}\label{eq7.22}
\end{equation}
where $c_{\ell'm'}$ are the coefficients of the spherical harmonics
expansion
\begin{equation}
\mathbf{e}^2 = \sum_{|m'|\le\ell'} \ c_{\ell'm'}\ Y_{\ell'm'}. \label{eq7.23}
\end{equation}
Summarizing (\ref{eq7.20})--(\ref{eq7.23}) we have the explicit
expression
\begin{equation}
\lambda''(0) = \lambda_\ell \sum\limits_{\underset{|m'| \le \ell'}{\ell'\neq\ell}} \
\frac{\lambda_{\ell'}}{\lambda_\ell-\lambda_{\ell'}} \
\big|c_{\ell'm'}\big|^2. \label{eq7.24}
\end{equation}
The coefficients $\big|c_{\ell'm'}\big|^2$ have been tabulated in \cite{FiedMisch}
 and are consistent with the detailed expansions
\begin{equation}
\mathbf{e}^2_{O(2)^-} = \frac{1}{ \sqrt{\pi}} \bigg(\frac12\ Y_{00} +
\frac{1}{\sqrt{5}}\ Y_{20}\bigg)\label{eq7.25a}
\end{equation}
for $\ell=1$, and 
\begin{equation}
\label{eq7.25b}
\begin{array}{lll}
\mathbf{e}^2_{O(2)^-} &=& \frac{1}{2\sqrt{\pi}} \left(Y_{00} +
\frac{4}{3\sqrt{5}}\ Y_{20} + \frac{6}{11}\ Y_{40} +
\frac{100}{33\sqrt{13}}\ Y_{60}\right);\\[2mm]
\mathbf{e}^2_{O^- }&=& \frac{1}{2\sqrt{\pi}} \left(Y_{00} - \frac{7}{11}\ Y_{40} +
\frac{\sqrt{70}}{11} \ {\rm Re} \ Y_{44} + \frac{10}{11\sqrt{13}}\
Y_{60} + \frac{10\sqrt{14}}{11\sqrt{13}} \ {\rm Re} \ Y_{64} \right)\\[2mm]
\mathbf{e}^2_{D^{d}_{6}} &=& \frac{1}{2\sqrt{\pi}} \left(-Y_{00} +
\frac{\sqrt{5}}{3}\ Y_{20} - \frac{3}{11} \ Y_{40} +
\frac{5}{33\sqrt{13}}\ Y_{60} + \frac{10\sqrt{7}}{\sqrt{429}} \ {\rm
 Re} \ Y_{66}\right)
\end{array}
\end{equation}
for $\ell=3$. Insertion into (\ref{eq7.24}) proves the remaining $\ell
= 3$ entries of Table \ref{tableofnonvanishing} and completes the proof of the lemma.

\hfill$\bowtie$

The bifurcating branches $(\lambda_K(s), v_K(s))$ of Lemma \ref{lm7.1}
are contained in any local center manifold graph $\Psi c$ at $\lambda
= \lambda_\ell$, $v =0$ as equilibria of the reduced ODE $\dot{v}_c =
\Phi(\lambda,v_c)$; see (\ref{eq7.5}). By $\Gamma$-equivariance, the bifurcating equilibria
$v_K(s)$ of isotropy $K$, for $s\neq 0$, generate group orbits $\Gamma
v_K(s) \cong \Gamma/K$ of equilibria which, likewise, appear as $C^k$
manifolds in Graph $\Psi c$. The linearized semiflow
\begin{equation}
\label{eq7.26}
\partial_t w = \Delta w +\lambda f'(v_K (s)) w
\end{equation}
of (\ref{6.4a}) at $v_K(s)$ therefore possesses a trivial eigenvalue 0 of geometric
multiplicity $\dim \Gamma/K$. This trivial part of the spectrum  is likewise contained in the center part  $\sigma_c$ which arises by
perturbation of the eigenvalue 0 of algebraic and geometric
multiplicity $\dim V_\ell = 2\ell+1$ of the linearized semiflow (\ref{eq7.26})
 at $\lambda = \lambda_\ell, \ v=0$ for $s=0$. By
standard linear perturbation theory this center part
$\sigma_c(B^\lambda(s))$ can be calculated from the linearization of
the reduced ODE $\dot{v}_c = \Phi(\lambda,v_c)$.

The remaining non-center part  $\sigma_h$ of the
linearization \ref{eq7.26} is strictly hyperbolic, and is inherited
from the strictly hyperbolic part at $s=0$. This decomposition leads
to the following proposition.

\setcounter{proposition}{1}
\begin{proposition}\label{prop7.2}
{~}

Let $(\lambda_K(s),v_K(s))$ bifurcate from $\lambda=\lambda_\ell, \
v=0$ as in Lemma \ref{lm7.1} above and consider small enough $\vert s\vert \neq 
0$.

Then the manifold $\Gamma v_K(s)$ of equilibria is normally hyperbolic for the semilinear semiflow (\ref{6.4a}) 
if, and only if, it is normally hyperbolic in the center
manifold. The strong unstable dimension $i=i(v_K (s))$  in $X = C^{\beta}(S^2)$, alias the Morse index,
then relates to the strong unstable dimension $i_c$ within the center
manifold as 
\begin{equation}
i= \ell^2 + i_c. \label{eq7.27}
\end{equation}
Normal hyperbolicity, with the same Morse index $i(v_K (s))$ as given by (\ref{eq7.27}), likewise holds true for the original quasilinear semiflow  (\ref{eq:v1.*} ), (\ref{v1.*in4} ),  (\ref{eq7.1} ) at the same equilibria.  
\end{proposition}

{\bf Proof:}
{~}
Morse equivalence lemma \ref{lm7.1} reduces the quasilinear case to the semilinear case. By standard center manifold reduction of the semilinear case, 
it only remains to prove that $\ell^2$ is the dimension of the
strong unstable part of the linearization \ref{eq7.26} at $s=0,\
\lambda=\lambda_\ell,\ v=0$. Since the strong unstable eigenspace of Re spec$<0$ is
$V_0 \oplus \ldots\oplus V_{\ell-1}$ with $\dim V_{\ell'} = 2\ell'+1$,
this latter claim is obvious since antiquity, and the proposition is proved.

\hfill$\bowtie$

We remark that the strong unstable dimension $i$ in (\ref{eq7.27})  is
also the codimension of the center-stable manifold of $(\lambda_K(s),
v_K(s))$. It coincides with the codimension of the stable set of the
group orbit $\Gamma v_K(s)$, which is foliated by the strong stable
manifolds over the base manifold $\Gamma/\Gamma_K$ of equilibria.
Although these claims which are standard for semilinear semiflows, also hold in the quasilinear semigroup setting, we do not pursue such details here.
\begin{theorem}\label{th7.3}
{~}

Consider the branches $(\lambda_K(s), v_K(s))$ of nontrivial
equilibrium solutions of (\ref{eq:v1.*}) which bifurcate at $\lambda_K(0) =
\lambda_\ell, \ v_K(0) = 0, \ \ell=1,2,3,4,$ and possess maximal
isotropy $K$.

Then the bifurcating group orbits $\Gamma v_K(s) \cong \Gamma/K$ of
equilibria are normally hyperbolic. Their strong unstable dimensions
are listed in Table \ref{tabletypes}. The local bifurcation diagrams are sketched
in Figure \ref{fig7.1}.
\end{theorem}

\psfrag{X}{$X$}
\psfrag{1}{$[1]$}
\psfrag{4}{$[4]$}
\psfrag{6}{$[6]$}
\psfrag{9}{$[9]$}
\psfrag{16}{$[16]$}
\psfrag{20}{$[20]$}
\psfrag{21}{$[21]$}
\psfrag{25}{$[25]$}
\psfrag{l=0}{$\ell=0$}
\psfrag{l=1}{$\ell=1$}
\psfrag{l=2}{$\ell=2$}
\psfrag{l=3}{$\ell=3$}
\psfrag{l=4}{$\ell=4$}
\psfrag{O(2)(2)}{$O(2)$, $[2]$}
\psfrag{O(2)+Z2}{$O(2)\oplus\mathbb{Z}^2$, $[5]$}
\psfrag{O-10}{$O^-$, $[10]$}
\psfrag{O(2)-12}{$O(2)^-$, $[12]$}
\psfrag{D6d13}{$D_6^d$, $[13]$}
\psfrag{O(2)+Zc2}{$O(2)\oplus\mathbb{Z}^c_2$, $[20]$}
\psfrag{lambda}{$\lambda$}
\psfrag{O+Zc2}{$O\oplus\mathbb{Z}^c_2$, $[18]$}
\setcounter{figure}{0} 
\begin{figure}[h]
\includegraphics[width=0.9\textwidth]{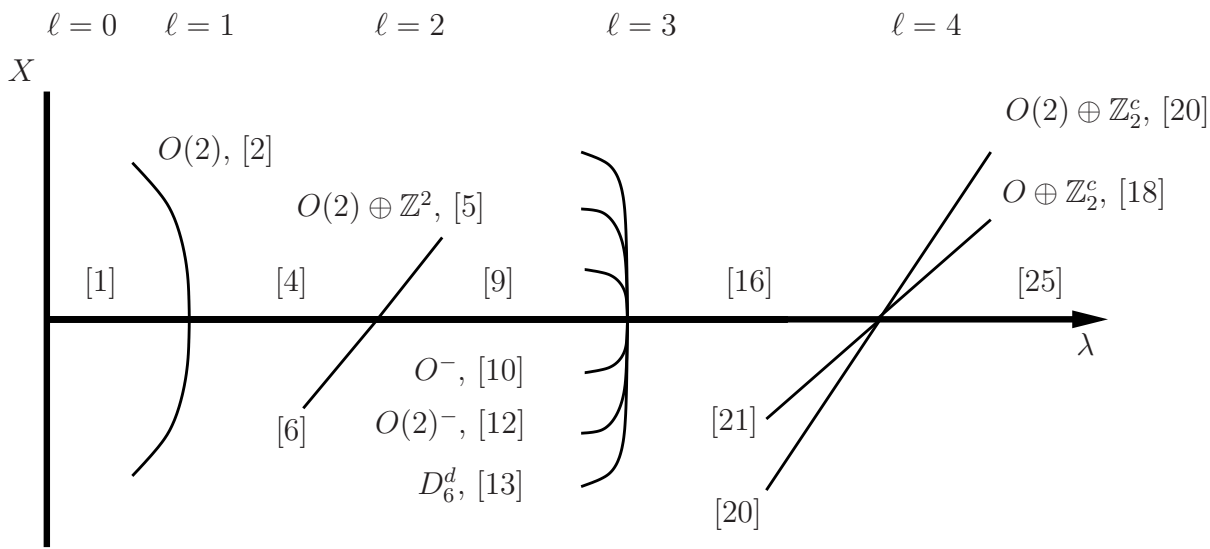}
\caption{\em{local transcritical and pitchfork bifurcation branches
   $(\lambda_K(s),v_K(s))$ from $\lambda_k(0) = \lambda_\ell =
   \ell(\ell+1), \ v_K(0) = 0$  with maximal isotropies $K$. Numbers
   in brackets indicate codimensions of their normally hyperbolic
   group orbits}.\label{fig7.1} }
\end{figure}

{\bf Proof:}
{~}

As we have seen in Proposition \ref{prop7.2}, normal hyperbolicity can
be checked within the center manifold. For $\ell\le 4$, as we consider
here, the local stability analysis of bifurcations with maximal
isotropy has been performed in \cite{LautHabil, CLM} in terms of lowest order
equivariant polynomials of quadratic and cubic order, respectively,
for even and odd $\ell$, and under certain not overly explicit
assumptions on genericity. Closer inspection shows that genericity,
in our cases, amounts to the bifurcation slopes $\lambda'_K(0)$ and
curvatures $\lambda''_K(0)$ to be nonzero and, for $\ell=3$, to avoid the
absence of the nontrivial cubic invariant.

Nonzero slopes and curvatures have been proved in Proposition
\ref{prop7.2}. To address cubic nondegeneracy for $\ell = 3$, we use \cite{CLM}, Table 6.2
 from which the bifurcation curvatures derive as
\begin{equation}\label{eq7.28}
\lambda'' _K(0) = \left\{
\begin{array}{lcl}
2(18\alpha-\beta )&  & K=O(2)^-\\[1.5mm]
8(45\alpha-2\beta) & \rm{for}& K=D^d_6\\[1.5mm]
16(10\alpha-\beta) && K=O^-.
\end{array}\right.
\end{equation}
Here the general cubic equivariant is written as $\alpha
\tilde{c}(v_c) + \beta|v_c|^2_2v_c$ in \cite{CLM}, with a suitable
equivariant vector polynomial $\tilde{c}(v_c)$ of the center manifold
variable $v_c \in V_\ell = \mathbb{R}^7$. 

The three bifurcation curvatures $\lambda''_K(0)$ of Lemma \ref{lm7.1} and
Table \ref{tableofnonvanishing} consistently overdetermine the two coefficients $\alpha,\beta$ as
\begin{equation}
\alpha= \frac{287}{1430\pi}, \qquad \beta = \frac{812}{143\pi}.\label{eq7.29}
\end{equation}
In particular $\alpha \neq 0$. By \cite{CLM}, Tables 6.2, 6.3, this
establishes normal hyperbolicity of all bifurcating branches
$(\lambda_K(s), v_K(s))$, for small $|s| \neq 0$, and the strong
unstable dimensions listed under $i_c$ in Table \ref{tabletypes}. These results
persist under perturbations by higher order terms; see \cite{CLM}, Theorem 5.2. 

For illustration consider the axisymmetric case $K= O(2)^-$ for
$\ell=3$. By Lemma \ref{lm7.1}, Table \ref{tableofnonvanishing}, we encounter a subcritical
pitchfork, $\lambda''_K(0) < 0$, which contributes a simple unstable
eigenvalue $\mu_0 > 0$ to the linearization at $(\lambda_K(s),v_K(s))$
in the center manifold. By \cite{CLM}, Table 6.2, the four remaining
directions transverse to the two-dimensional group orbit
\begin{equation}
\Gamma v(s) \cong \Gamma/\Gamma_{v(s)} = O(3)/O(2)^- = SO(3)/SO(2) =
S^2 \label{eq7.30}
\end{equation} 
in $V_\ell \cong \mathbb{R}^7$ produce two further unstable eigenvalues $\mu_1, \mu_2 >0$, for cubic coefficients $\alpha\neq 0$. This proves $i_c=3$ and, by Proposition \ref{prop7.2}
\begin{equation}
\label{7.31}
i=i_c+\ell^2=12
\end{equation}
as claimed in Table \ref{tabletypes}.

The remaining entries in Table \ref{tabletypes} are obtained similarly, with group orbits $$\Gamma v\cong O(3)/ \left(  O(2) \oplus \mathbb{Z}_2^c\right) = \mathbb{R}P^2$$ in the axis symmetric cases $\ell=2,4$ and with 3-dimensional homogeneous spaces $\Gamma v$ for the discrete subgroups $K$.
 
The local bifurcation diagrams of Figure \ref{fig7.1}  combine the
information of Table \ref{tableofnonvanishing} on the bifurcation types and directions with
the unstable dimensions $i$ of Table \ref{tabletypes}. This proves the theorem.

\hfill$\bowtie$

We can now combine the results of the previous chapters to address normal hyperbolicity and the strong stable manifolds of the symmetry breaking equilibrium branches $\left(\lambda(s), v_K(s)\right) $ of \ref{eq:v1.*}, which arise in Lemma \ref{lm7.1}, Theorem \ref{th7.3} and Table \ref{tableofnonvanishing}. \\
\begin{table}[h!]
\begin{center}
\begin{tabular}{||c|c|c|c|c|c||} \hline\hline
$\ell$ & isotropy $K$ & type & direction &$i_c$ & $i$  \\ 
\hline\hline 
$1$ &$  O(2)^- $ & pitchfork & sub & 1 & 2  \\
 \hline \hline
 & && sub & 2 & 6 \\
\cline{4-6}
2&$O(2)\oplus \mathbb{Z}_2^c$&transcritical&super&1 & 5\\
\hline \hline
 & $O(2)^-$& pitchfork &sub& 3  &12   \\
 \cline{2-6} 
 $3$&$O^-$&pitchfork&sub&1&10\\
 \cline{2-6}
&  $D_6^d$& pitchfork &sub&  4& 13 \\
\hline \hline
 & &&  sub&  4 & 20\\
\cline{4-6}
$4$&$O(2)\oplus \mathbb{Z}_c^2$&transcritical&super&4&20\\
 \cline{2-6} 
 &&&sub&5 &21\\
 \cline{4-6}
 &$O^-\oplus \mathbb{Z}_c^2$&transcritical&super&2&18\\
\hline \hline

\end{tabular}
\end{center}
\caption{For bifurcations $\left(\lambda(s), v_K(s)\right) $ from $\lambda_K(0)=\lambda_{\ell}=\ell (\ell +1)$, $v_K(0)=0$ with maximal isotropies $K$ and normally hyperbolic group orbits,  we indicate the bifurcation type, direction, strong unstable dimension $i_c$ within the semilinear center  manifold and strong unstable quasilinear dimension $i$ in $X=C^{\beta}(S^2)$, alias the codimensions of the strong stable manifolds of the group orbits $\Gamma v_K(0)$ for $s\neq 0$. \label{tabletypes}
} 
\end{table}

\section{Conclusions and outlook}
\label{concl}

The existence of a manifold $\Gamma v_K \cong \Gamma/K$ of nontrivial equilibria of the self-similarly rescaled equation (\ref{eq:v1.*}), 
as well as their strong stable manifold $W^{ss}(\Gamma v_K)$ has been studied above in the context of equivariant bifurcation theory. 
The normal hyperbolicity of the manifolds $\Gamma v_K$ has been established in Theorem \ref{th7.3}. By Theorem \ref{th5.1}, 
we conclude  the existence of solutions $v$ of the self-similarly rescaled equation (\ref{eq:v1.*}) which converge to a branching 
equilibria $v(t, \cdot)\rightarrow v_K$ as $t \rightarrow +\infty$. 
 In the original equation \ref{eq:parscal} this behavior corresponds to an asymptotically self-similar blow up as $r\rightarrow r_1=1$ with 
\begin{equation} \label{uvchange}
 u(r,p)  = \left( \frac{\lambda}{2}\left(1/r -1\right)\right)^{-\frac12}\left( v(-\frac{1}{2}\log(1-r), p) + 1\right), \qquad p\in S^2,
 \end{equation}
 with isotropic scalar curvature $R(r)=\frac{\lambda +2}{r^2}$.
\begin{remark} \label{heteroclinics} (heteroclinic orbits)\\
Equivariant bifurcation theory does not only provide a bifurcation diagram of branches of 
equilibria of the self-similarly rescaled equation (\ref{eq:v1.*}) as given in Figure \ref{fig7.1}. 
As Theorem \ref{centermnf} provides a local center manifold, we deduce the existence of heteroclinic orbits connecting  
those equilbria. Following  Lauterbach et al. \cite{LautHabil, CLM},   recalling the fact that our equation satisfies the 
genericity condition required there as proved in the previous section (see proof of Theorem (\ref{th7.3})), 
and given the signs of the first nonvanishing derivatives in Table \ref{tableofnonvanishing}, we conclude the existence of heteroclinic 
connections between the following branches of the bifurcation diagram \ref{fig7.1} as summarized in the following Table \ref{tableheteroclinics}, 
where arrows $\rightarrow$ denote the existence of heteroclinics orbits between the branches of equilibria.\\
\end{remark}
\begin{table}[h!]
 \begin{center}
 \begin{tabular}{||c|c|c|c||}
 \hline \hline
 $\ell$ & type & $\lambda< \lambda_{\ell}$ & $\lambda > \lambda_{\ell}$ \\
 \hline \hline
 1 & pitchfork & $O(2)\rightarrow 0$ & - \\
 \hline
 2 & transcritical &$O(2)\oplus \mathbb{Z}_c^2\rightarrow 0$&$0 \rightarrow O(2)\oplus \mathbb{Z}_c^2$\\
 \hline
 3 & pitchfork & $D_6^d \rightarrow  O(2)^- \rightarrow O^- \rightarrow 0$& - \\
 \hline
 4 & transcritical &$O(2)\oplus \mathbb{Z}_c^2\leftarrow O\oplus \mathbb{Z}_c^2\rightarrow 0$& $0\rightarrow O\oplus \mathbb{Z}_c^2\leftarrow O(2)\oplus \mathbb{Z}_c^2$\\
 \hline\hline
 \end{tabular}
 \end{center}
 \caption{Heteroclinic orbits between bifurcating branches. \label{tableheteroclinics} }
 \end{table}

\begin{figure}[h]
\includegraphics[width=0.9\textwidth]{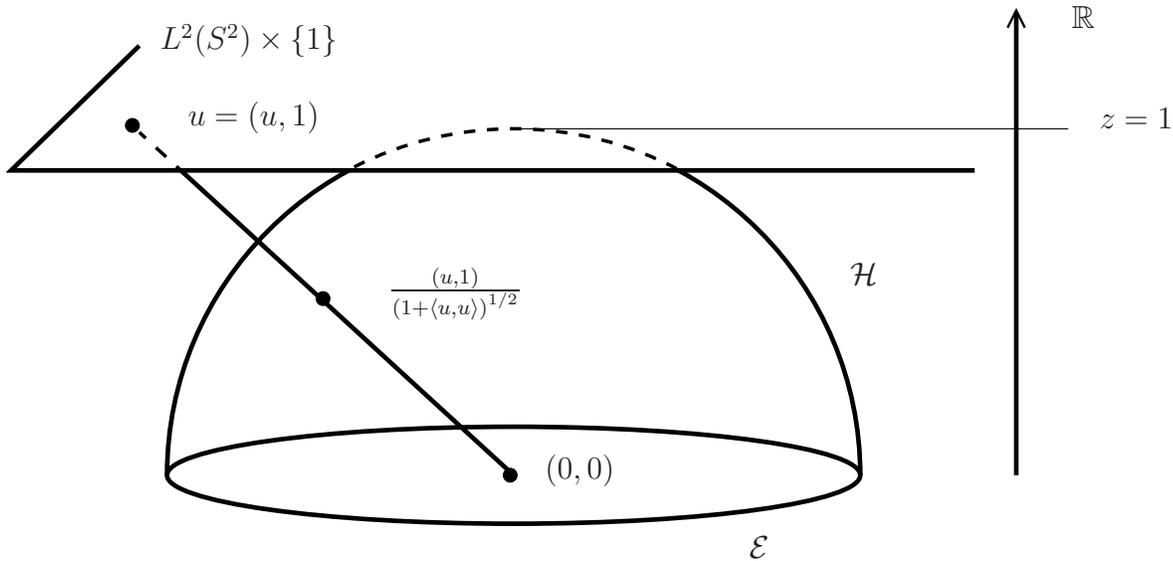}
\caption{\label{compactification}  {\em  Poincar\'e compactification: each element $u\in L^2(S^2)$ is identified with $(u,1)\in L^2(S^2)\times \{ 1 \} $ and projected centrally on the unit upper  Hemisphere $\mathcal{H}$. Its equator $\mathcal{E}$ is the sphere at infinity.    }}
\end{figure}

The self-similar rescaling by the ODE blow up rate $\left( 1/r-1\right)^{-\frac12}$ can be also  seen as a Poincar\'e ``compactification" that allows us to describe the blow up behavior as the convergence to an equilibrium in the sphere at infinity. 
To be more precise, the Poincar\'e compactification projects centrally the ambient Hilbert space $L^2(S^2)$ into the  unit hemisphere  $\mathcal{H}$ of $L^2(S^2)\times \mathbb{R}$ (which is a ball of the same dimension as $L^2(S^2)$ itself), in such a way that infinity is  projected on its boundary,  the equator $\mathcal{E}$, alias the sphere at infinity, whose elements will be denoted by $\chi\in L^2(S^2)$, $\Vert \chi \Vert =1$ in the following. See Figure \ref{compactification} for a sketch of the Poincar\'e compactification and \cite{Hell} for details.
\psfrag{XX}{$L^2(S^2)\times \{1 \}$}
\psfrag{v}{$u=(u,1)$}
\psfrag{pv}{$\frac{(u,1)}{\left( 1+ \langle u,u \rangle   \right)^{1/2}}$} 
\psfrag{R}{$\mathbb{R}$}
\psfrag{z=1}{$z=1$}
\psfrag{000}{$(0,0)$}
\psfrag{Einf}{$\mathcal{E}$}
\psfrag{H}{$\mathcal{H}$}

\begin{figure}[!h]
\begin{center}
\includegraphics[width=0.7\textwidth]{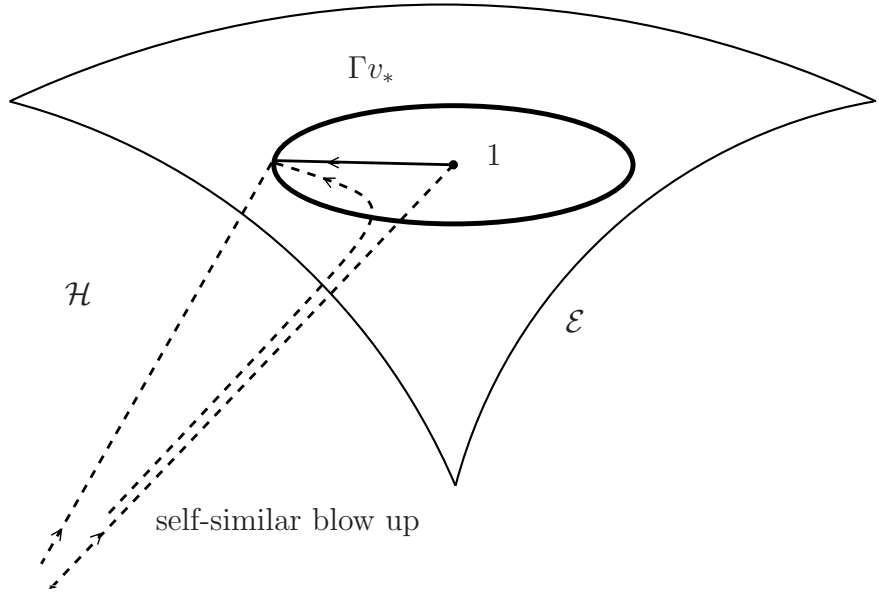}
\end{center}
\caption{\label{urbietorbi}  {\em This Picture caricatures the infinite dimensional Poincar\'e hemisphere $\mathcal{H}$ of $L^2(S^2)$ in 3-d. 
The single point denoted by 1 is the isotropic equilibrium at infinity. The ellipse denoted by $\Gamma v_*$ stands for a 
manifold of anisotropic equilibria. The straight dashed lines coming from the interior of the Poincar\'e hemisphere are self-similar
blow up solutions with profile 1 and $v_*$ respectively. The sphere $\mathcal{E}$ at infinity contains a  heteroclinic  $1 \rightarrow v_*$ that influences the blow up behavior of non self-similar blow up solutions.  }}
\end{figure}

 This projection is called compactification for historical reasons - it has been developed by Poincar\'e for ODE's - 
but in the case of function spaces such as $L^2(S^2)$, the resulting Poincar\'e hemisphere $\mathcal{H}$ is not compact. 
Following \cite{Hell}, the Poincar\'e compactification induces a compactified equation on $\mathcal{H}$ as well as an equation on the 
sphere $\mathcal{E}$ at infinity, which is the one we are interested in.\\
We begin from the equation (\ref{eq:parscal}), which can be simplified  in our simple settings  by using $u\mapsto r^{-1/2}u$ to
\begin{equation} \label{simpleeq}
\partial_r u= u^2\Delta u + \frac{\lambda}{2} u^3. 
\end{equation} 
See \cite{S} for details. 
Then the Poincar\'e compactification induces after normalization the following equation on the sphere at infinity:
\begin{equation}
 \label{compeqE} \partial_{\tau} \chi =  \chi^2\Delta \chi + \frac{\lambda}{2}\chi^3  - \left\langle  \chi^2\Delta \chi + \frac{\lambda}{2} \chi^3, \chi \right\rangle  \chi , \qquad \Vert \chi\Vert=1, \qquad \chi \in \mathcal{E}, 
 \end{equation}
 where $\left\langle.,. \right\rangle $ is the scalar product on $L^2$ and the time variable $\tau$ corresponds to a normalization 
preventing the trajectories from hitting   the sphere at infinity in finite time. 
The previous equation (\ref{compeqE}) is  not a full PDE because it contains nonlocal terms coming from the scalar product. 
Now we can interpret the self-similarly rescaled equation (\ref{eq:nueq}) ( which differs from (\ref{eq:v1.*}) only by a 
multiplication by a scalar and a translation of the origin) as an equation on the sphere at infinity. 
Let $\nu$ be a $L^2(S^2)$-bounded solution of equation (\ref{eq:nueq})
\[\partial_t  \nu= \nu^2 \Delta \nu  + \frac{\lambda}{2}\nu^3-\nu.\]
Then $\chi= \frac{\nu(t, p)}{\Vert \nu(t,. ) \Vert}\in \mathcal{E}$ solves the equation on the sphere at infinity (\ref{compeqE}) 
up to a factor that can be gotten rid of by another change of time variable. More precisely, 
computed in the time variable $t$ defined by $r=1-\exp(-2t)$, the above $\chi$ satisfies the equation 
\begin{equation}
\label{insphereinfinity}
\partial_t \chi = \langle \nu(t,.),   \nu(t,.) \rangle \left( \chi^2\Delta \chi + \chi^3  - \left\langle  \chi^2\Delta \chi 
+ \chi^3, \chi \right\rangle  \chi \right)
\end{equation}
The heteroclinic connections  discussed previously are objects in the sphere $\mathcal{E}$ at infinity, 
and influence the blow up  behavior of our original Problem as depicted in Figure \ref{urbietorbi}, 
which is a 3-d caricature of this infinite dimensional hemisphere (ball) $\mathcal{H}$. 
Here we see that  trajectories with finite initial conditions in the interior of  the Poincar\'e hemisphere $\mathcal{H}$  
are caught by the heteroclinic connection between the trivial equilibrium at infinity and  a manifold $\Gamma v_*$ of 
anisotropic equilibria at infinity,  and hence explode anisotropically. 

\psfrag{selfsim}{self-similar blow up} 
\psfrag{gammma}{$\Gamma v_*$}
\psfrag{111}{1}

Now that we have sketched the dynamics at infinity of the original Equation (\ref{eq:parscal}), 
let us interpret its geometrical  meaning  for the manifold $M$ at its upper boundary $\{ r_1\} \times S^2$. 
In fact, the blow-up behavior reflects a degeneracy of the area-radius $r$, but not a degeneracy of the metric $g$ 
of the manifold $M$ at its upper boundary $\{ r_1\} \times S^2$. 
Let us recall that the metric on the manifold $M$ is given by 
\begin{equation}
\label{metric}
g= u^2 dr^2 + r^2 \omega,
\end{equation} 
 where $\omega$ is the standard metric of the round sphere. Now we introduce a coordinates $s$ proportional to the  geodesic length as
 \begin{equation}
 \label{scoord}
 s:=\int_{r_0}^r \left( 1/\bar{r}-1\right)^{1/2} d\bar{r}.
 \end{equation}
 Hence 
 \begin{equation}
 dr = (1/r -1)^{\frac{1}{2}} ds,
 \end{equation}
 and the metric belonging to the  (asymptotically) self-similar blow-up $u=\left( \frac{\lambda}{2}\left(1/r -1  \right)\right)^{-\frac{1}{2}} v$ 
reads in these coordinates:
 \begin{equation}
 \label{nonsingmetric}
 \begin{array}{lll}
 g &  =  & \left( \frac{\lambda}{2}\left(1/r -1  \right)\right)^{-1}   v^2  (1/r -1)ds^2 + r^2(s) g^{\Sigma}\\
 & & \\
  & =&\frac{\lambda}{2} v^2ds^2+ r^2(s) g^{\Sigma}
  \end{array}
 \end{equation}
 In this form, the metric is continuous and nondegenerate  up to  and including the  upper boundary of the 
manifold $M$: $g\in C^{\infty}(M)\cap C^0(\bar{M})$, no matter which profile for $v$ is used. It may be the trivial 
equilibrium of the self-similarly rescaled Equation (\ref{eq:v1.*}), one of the equilibria or heteroclinics we exhibited 
nearby via equivariant bifurcation theory, or any solution converging to any other  equilibrium $v$  of the self-similarly 
rescaled  Equation (\ref{eq:v1.*}).
 
The forward time evolution of a heteroclinic connection between two equilibria $v_*^1$, $v_*^2$ of Equation (\ref{eq:v1.*}) 
has just been discussed. Let us now discuss the backward time evolution 
$t\rightarrow - \infty$ alias $r\rightarrow - \infty$ through $r=1- \exp(-2t)$ of a heteroclinic trajectory $v(t)$ from 
equilibrium $v_*^1$ to equilibrium $v_*^2$. Replacing the latter change of time/radius variable by $\tilde{r}=1-\exp(-2\theta)\exp(-2t)$ 
corresponds to a time shift $t\mapsto t+\theta$. Let $r_0$ be the initial radius of our construction. 
The first change of radius/time variable puts $r_0$ at time $t_0=\ln \left(   (1-r_0)^{-1/2} \right)$, 
while the second puts $r_0$ at $t_0+\theta$. In the second case, $v(t_0+\theta) $ determines the shape of the metric 
at initial radius $r_0$: $v(t_0+\theta)$ comes arbitrarily near equilibrium $v_*^1$ as $\theta$ approaches $-\infty$. 
Once $\theta$ is fixed, the equilibrium $v_*^2$ at the other end of the heteroclinic gives the shape of the metric at the blow up 
radius $r_1=1 \Rightarrow t, t+\theta \rightarrow \infty$ for any fixed $\theta$.

Let us now discuss  the question of the center  of $M$, that is, the limit where the lower boundary disappears as  $r_0 \rightarrow 0$. 
Because we chose to study bifurcation with respect to the parameter $\lambda$, this imposes the prescribed scalar curvature to be 
$R(r)=\frac{\lambda +2}{r^2}$. Therefore, the scalar curvature $R$, and with it at least one sectional curvature, concentrates  at $r=0$: 
we cannot expect $M$ to be regular at the center.\\
This is seen again by considering the coordinate $s$ proportional to the geodesic length, and for simplicity the isotropic self-similar 
blow-up profile $v\equiv 1$, with blow up radius $r_1=1$. Near $r=0$, the radial component of the metric can be approximated by  
$\frac{r}{1-r}\approx r$, so that the metric and the geodesic length are given by
\begin{equation}
\begin{array}{lll}
g & = &  rdr^2 + b r^2\omega \\
s &=&  \int_0^r \sqrt{\bar{r}} d\bar{r} = \frac{2}{3} r^{\frac{3}{2}} \\
g &=& ds^2 + b s^{\frac{4}{3}} \omega
\end{array}
\end{equation} 
The singularity manifests itself in the term $s^{\frac{4}{3}}$: the discrepancy between the exponent $\frac{4}{3}$ and $2$ 
excludes regularity of the center. \\
Our  motivation was to construct $M$ as  a piece of initial conditions for the Einstein equations. In view of this,  
curvature singularity at the center may not be  desirable. To avoid such  a situation, one may want to fix a $r_0>0$ and  
use a different  method for prescribing  the scalar curvature for $r\leq r_0$. 
%

 Even if the  isotropic scalar  curvature $R(r)=\frac{\lambda+2}{r^2}$, which we have prescribed in this work for 
simplicity in the arguments involving bifurcation theory  is rather primitive,    we are still able to show that a rich 
anisotropic and only asymptotically self-similar behavior develops.  The isotropy of the metric component in radial direction 
can break in several ways as the boundary $\{r_1\}\times S^2$ of the manifold $M$ is approached, depending on the choice of the 
bifurcation parameter $\lambda$.   This method may be applied to more realistic situations where $\omega$ and $r^2 R $ 
do depend on $r$. In this case, 
  it may be possible to apply the equivariant bifurcation theory arguments  on a thin collar $[r_0, r_1)\times S^2$ region, 
and complete the manifold by prescribing a scalar curvature leading to a regular center in the center region $[0, r_0]\times S^2$.\\ 

Let us briefly address the behavior of the metric $g$ at $r=0$, according to our construction. 
The parabolic scalar curvature equation \ref{eq:parscal} is equivariant, separately, 
with respect to each of the involutions $r \leftrightarrow -r \in [-1, 1]$ and $ \varphi \leftrightarrow - \varphi \in \Sigma=S^2$. 
Our transformation $u\mapsto \nu$ to the self-similar rescaled equation (\ref{eq:nu,r}) preserves both involutions if we define 
\begin{equation}
\label{rescrneg}
\nu = \sqrt{ 1/ \vert r \vert -1} \quad u,
\end{equation} 
to accomodate negative $r$. A simple reflection 
\begin{equation} \label{reflection}
\nu(-r,p):=\nu(r,p),
\end{equation}
for $r\geq 0$ then provides a consistent extension to negative radial arguments. In terms of smooth elements $v$ 
in a strong stable manifold $W^{ss}$ of (\ref{eq:v1.*}) we get
\begin{equation}
\label{unu}
u=(1/\vert  r \vert -1)^{-1/2} \nu = \vert r \vert ^{1/2} ( 1+\frac12 \vert r \vert + \ldots )\nu.
\end{equation}
With $\nu=\left(\lambda/2\right)^{-1/2} (v+1)$, and $v$ regular in $t=r+\dots$ at $r=0$ we obtain an expansion
\begin{equation}
\label{expansion}
u^2=\vert r \vert (1+\vert r \vert +\ldots)\frac{2}{\lambda} (1+v)^2.
\end{equation}
This identifies the metric
\begin{equation} \label{met}
g=u^2 dr^2 + r^2 \omega
\end{equation}
to possess a quadratic cusp at $r=0$. Nonetheless, Definition (\ref{reflection}) provides a valid extension of $u$ beyond $r=0$ into negative $r$. 

Suppose we then insist, in addition, that $(r,p)$ and $(-r,-p)$ denote the same point on the 3-manifold $M=[0,r_1]\times S^2$ as 
common in polar coordinates at the point $r=0$ where the  cusp singularity resides. Then we have to require   antipodal symmetry 
\begin{equation}
\label{antipodal}
\nu(r,-p)=\nu(-r, p)= \nu(r,p)
\end{equation}
 for $\nu$ and consequently for $v$. This requirement amounts to isotropy of $v$ under $-\rm{id} \in O(3)$, in the language of Section \ref{sym}: 
\begin{equation} \label{incl-id}
\Gamma_v \geq \mathbb{Z}_2^c = \langle -\rm{id} \rangle. 
\end{equation}
By our $O(3)$ bifurcation analysis of Section 3, this requirement is automatically satisfied for all symmetry breaking branches 
$(\lambda(s), v(s))$ which emanate at $\lambda(0)=\lambda_{\ell}= \ell (\ell+1)$ with even $\ell=0,2,4,\ldots$. But (\ref{antipodal}),(\ref{incl-id}) 
is just as automatically violated at odd $\ell$. From this point of view, even $\ell$ are preferable because they do not need any further 
treatments for $0\leq r \leq r_1$ and because the simple reflection (\ref{reflection}) defines a consistent, and symmetric 
extension to $-r_1 \leq r \leq 0$.  

As noted previously, the initial data for the Einstein equations have to be completed by a) a regular center, and b) a region matching 
the blow up boundary $\{r_1\}\times S^2$ with  the apparent horizon of the black hole. The latter can be constructed with anisotropy, 
for more details see \cite{Smith1, Smith2}. Both  sides of the matching region,   $\{r_1\}\times S^2$ and the apparent horizon are critical 
points of the area functional, the latter being a minimum of Morse index 0, and the former a critical point of index greater than 1. 
To see this, let us notice that at blow up radius $r=1$ , the mean curvature $H=\bar{H}/ru \propto (1-r)^{1/2}/r^{3/2}$ vanishes, 
and this implies that the sphere $S^2\times {1}$ is a critical point of the area functional  $\mathcal{A}$, i.e. a so called minimal surface.
If we denote by $F_{\tau}:\Sigma =S^2\rightarrow M$ a smooth family of spheres parametrized by a parameter $\tau$, such that $\partial_{\tau} F $ 
is a multiple of the outer normal $N$ to $\Sigma_{\tau}=F_{\tau}(S^2)$ with factor $\varphi$, $\partial_{\tau} F =\varphi N$, then  (see \cite{chav})
\begin{equation}
\label{first-variation-area}
\frac{d\mathcal{A}}{d\tau}= - \int_{S^2} \varphi H d\mathcal{A}.
\end{equation}
To analyze the type of critical point we are facing, we need the second variation of the area functional at $\tau=0$, 
which reads in our simple settings
\begin{equation}
\label{second-variation-area}
\left. \frac{d^2\mathcal{A}}{d\tau^2} \right\vert_{\tau=0}= - \int_{S^2} \vert \nabla \varphi \vert^2 + \varphi^2 \left( 1-R/2 \right).
\end{equation}
The self adjoint operator associated to this functional is 
$\Delta - \left( 1- R/2\right)= \Delta + \lambda/2$, 
and the stability index $i_{\mathcal{A}}$ of a critical point of $\mathcal{A}$ is determined  by the number of positive 
eigenvalues (counted with their multiplicity) of this operator. Hence we get for $\lambda\in (2\lambda_{\ell}, 2\lambda_{\ell+1})$ 
that the stability index of the minimal surface at blow up radius  is $i_{\mathcal{A}}= 2\ell(\ell+1)+1$

Once the completions a) and b) have been performed,
initial conditions with anisotropic metrics will be provided for the Einstein equations. 
We will then be in a position to ask about their time evolution and study, in particular, how the anisotropy affects  the long  time behavior. 

\vspace*{3cm}
\bibliographystyle{amsplain}

\end{document}